\documentclass[11pt]{article}

\usepackage{amssymb}
\usepackage{amscd}
\usepackage{amsmath}
\usepackage{amsthm}

\renewcommand{\a}{{\alpha}}
\renewcommand{\b}{{\beta}}
\newcommand{\g}{{\gamma}}
\renewcommand{\d}{{\delta}}
\newcommand{\e}{{\varepsilon}}

\renewcommand{\l}{{\lambda}}

\newcommand{\f}{{\varphi}}
\renewcommand{\o}{{\omega}}

\renewcommand{\t}{{\theta}}
\newcommand{\s}{{\sigma}}
\renewcommand{\r}{{\rho}}

\newcommand{\D}{{\Delta}}
\newcommand{\G}{{\Gamma}}
\renewcommand{\L}{{\Lambda}}

\newcommand{\calS}{{\mathcal S}}

\DeclareMathOperator{\iindex}{index}
\DeclareMathOperator{\sign}{sign}

\DeclareMathOperator{\rank}{rank}

\DeclareMathOperator{\Dom}{Dom}
\DeclareMathOperator{\Ker}{Ker}
\DeclareMathOperator{\Tr}{Tr}
\DeclareMathOperator*{\slim}{s-lim}
\DeclareMathOperator{\spflow}{sf}

\renewcommand\Im{\hbox{{\rm Im}}\,}
\renewcommand\Re{\hbox{{\rm Re}}\,}

\newcommand{\abs}[1]{\left\lvert#1\right\rvert}
\newcommand{\norm}[1]{\left\lVert#1\right\rVert}
\newcommand{\snorm}[2]{\left\lVert#1\right\rVert_{\SS_{#2}}}
\newcommand{\br}[1]{\left(#1\right)}

\newcommand{\R}{{\mathbb R}}
\newcommand{\N}{{\mathbb N}}
\newcommand{\Z}{{\mathbb Z}}
\newcommand{\C}{{\mathbb C}}
\newcommand{\T}{{\mathbb T}}

\renewcommand{\H}{{\mathcal H}}
\newcommand{\K}{{\mathcal K}}
\newcommand{\calB}{{\mathcal B}}

\renewcommand{\SS}{{\mathfrak S}}

\numberwithin{equation}{section}

\theoremstyle{plain}
\newtheorem{theorem}{\bf Theorem}[section]
\newtheorem{lemma}[theorem]{\bf Lemma}
\newtheorem{proposition}[theorem]{\bf Proposition}
\newtheorem{condition}[theorem]{\bf Assumption}

\theoremstyle{definition}
\newtheorem{definition}[theorem]{\bf Definition}

\theoremstyle{remark}
\newtheorem*{remark}{\it Remark}

\newcommand{\thm}[1]{{Theorem \ref{t.#1}}}
\newcommand{\lma}[1]{{Lemma \ref{l.#1}}}
\newcommand{\prp}[1]{{Proposition \ref{p.#1}}}

\newcommand{\cnd}[1]{{Assumption \ref{condition.#1}}}
\newcommand{\dfn}[1]{{Definition \ref{d.#1}}}

\renewcommand{\qed}{\vrule height7pt width5pt depth0pt}

\newcommand{\xtp}{\widetilde X_p}
\newcommand{\rtp}{\widetilde \rho_p}
\newcommand{\wt}{\widetilde}
\newcommand{\Zo}{\Z_0}
\newcommand{\limn}{\lim_{n\to\infty}}
\DeclareMathOperator*{\xplim}{X_p-lim}

\DeclareMathOperator*{\xinflim}{X_\infty-lim}
\DeclareMathOperator*{\splim}{\SS_p-lim}
\DeclareMathOperator*{\sqlim}{\SS_q-lim}
\DeclareMathOperator*{\blim}{n-lim}
\DeclareMathOperator{\ext}{ext}

\newcommand{\J}{J^{-1}}

\title{The spectral shift function and the invariance principle}

\author{Alexander Pushnitski\thanks{
Department of Mathematical Sciences,
Loughborough University,
Loughborough, Leicestershire, LE11 3TU,
United Kingdom.
E-mail: a.b.pushnitski@lboro.ac.uk
}}

\date{November 1999}

\begin{document}

\maketitle

\begin{abstract}
The new representation formula for the spectral shift
function due to F.~Gesztesy and K.~A.~Makarov is considered. 
This formula is extended to the case of relatively trace
class perturbations.
\end{abstract}

%*************************************************************
\section{Introduction}
\label{section.a}
%************************************************************

{\bf 1.} 
First we briefly remind the definition of the spectral shift function 
(SSF). For the details and references to the literature, 
see \cite{BYa,Yafaev}.

Let $H_0$ and $H$ be self-adjoint operators in a Hilbert space 
$\H$, and let their difference belong to the trace class:
\begin{equation}
H-H_0\in\SS_1.
\label{a.1}
%................................................................(a.1)
\end{equation}
Then there exists a unique function  $\xi(\cdot; H,H_0)\in L_1(\R)$,
such that the following trace formula holds:
\begin{equation}
\Tr(\f(H)-\f(H_0))=\int_{-\infty}^\infty
\f'(\l)\xi(\l;H,H_0)d\l,
\quad\forall \f\in C_0^\infty(\R).
\label{a.2}
%................................................................(a.2)
\end{equation}
The function $\xi$ is called the SSF for the pair $H_0$, $H$. 

Let $\D_{H/H_0}(z)=\det((H-zI)(H_0-zI)^{-1})$, $\Im z>0$, 
be the perturbation determinant of the pair $H_0$, $H$. 
The following Krein's formula expresses the SSF 
in terms of the perturbation determinant:
\begin{equation}
\xi(\l;H,H_0)=\frac1\pi\lim_{y\to+0}\arg\D_{H/H_0}(\l+iy),
\label{a.3}
%................................................................(a.3)
\end{equation}
where the branch of the argument is fixed by the condition 
\begin{equation}
\lim_{y\to+\infty}\arg\D_{H/H_0}(\l+iy)=0.
\label{a.4}
%................................................................(a.4)
\end{equation}
The Birman-Krein formula relates the SSF to the scattering matrix 
$\calS(\l;H,H_0)$ for the pair $H_0$, $H$
(for the  definition of the scattering matrix, 
see, e.g., \cite{Yafaev}):
\begin{equation}
\det\calS(\l;H,H_0)=\exp(-2\pi i\xi(\l;H,H_0)),
\label{a.5}
%................................................................(a.5)
\end{equation}
for a.e. $\l$ on the absolutely continuous spectrum of $H_0$.

{\bf 2.} 
In \cite{GM}, a new representation for the SSF has been found.
In order to write down this representation, let us present
the perturbation $V:=H-H_0$ in the factorised form $V=G^* JG$, where
$G$ is a Hilbert-Schmidt operator, and $J=J^*=J^{-1}=\sign V$.
Further, denote
\begin{equation}
\begin{split}
A(\l+i0)&:=\lim_{y\to0+}\Re\bigl(G(H_0-(\l+iy)I)^{-1}G^*\bigr),\\
B(\l+i0)&:=\lim_{y\to0+}\Im\bigl(G(H_0-(\l+iy)I)^{-1}G^*\bigr).
\end{split}
\label{a.6}
%................................................................(a.6)
\end{equation}
Note that the limits in \eqref{a.6} exist for a.e. $\l\in\R$ 
in the operator norm
(and even in the norm of the Schatten-von Neumann ideal 
$\SS_p$ for any $p>1$ --- see \cite{BE,Naboko}).

The representation of \cite[Theorem 5.4]{GM} reads as follows:
\begin{equation}
\xi(\l;H,H_0)=\frac1\pi\int_{-\infty}^\infty\frac{dt}{1+t^2}
\iindex\bigl(E_{J+A(\l+i0)+tB(\l+i0)}((-\infty,0)),
E_J((-\infty,0))\bigr),
\quad
\text{ a.e. }\l\in\R.
\label{a.7}
%................................................................(a.7)
\end{equation}
Here $E_M(\cdot)$ stands for the spectral projection of a self-adjoint 
operator $M$, and $\iindex(\cdot,\cdot)$ denotes the index of a 
Fredholm pair of projections (see \eqref{b.1} below).
In the special case of perturbations of a definite sign
(where $J=\pm I$) the formula \eqref{a.7} 
was originally found in \cite{Push}.

{\bf 3.}
In applications, the assumption \eqref{a.1} 
becomes too restrictive.
Instead of \eqref{a.1}, it is usually possible to check that 
\begin{equation}
f(H)-f(H_0)\in\SS_1,
\label{a.8}
%................................................................(a.8)
\end{equation}
where $f:\s(H_0)\cup\s(H)\to\R$ is a monotone smooth enough function. 
In what follows, we for simplicity assume that $f$ is non-decreasing
(otherwise one can replace $f$ by $-f$). 

Under the assumption 
\eqref{a.8}, the SSF for the pair $f(H_0)$, $f(H)$
exists and the corresponding trace formula  is valid. 
The change of variables $\l\mapsto f(\l)$ leads to
the trace formula \eqref{a.2} for the pair $H_0$, $H$ with 
\begin{equation}
\xi(\l;H,H_0)=\xi(f(\l);f(H),f(H_0)).
\label{a.9}
%................................................................(a.9)
\end{equation}
Usually formula \eqref{a.9} is treated as the definition of the SSF 
$\xi(\cdot;H,H_0)$ under the assumption \eqref{a.8}.
Further details can be found in \cite[\S 8.11]{Yafaev}.
For the function $f$, one often takes 
$f(\l)=(\l-\l_0)^{-m}$ or $f(\l)=e^{-a\l}$.

{\bf 4.}
For the case of perturbations $V$ of a definite sign
and semibounded from below operators $H_0$, $H$, 
formula 
\eqref{a.7} has been extended (in \cite[Theorem 1.2]{Push})
to the case when the inclusion \eqref{a.8}
(but not necessarily \eqref{a.1}) holds true
with $f(\l)=(\l-\l_0)^{-m}$.
This extension has proved to be useful in applications 
to differential operators (see \cite{Push2}).

The aim of this paper is to prove a similar result 
without the assumption on the sign of the perturbation.
Below we briefly describe our main result; 
for a precise statement, see \thm{f.5b}.

Let $H_0$ be a self-adjoint operator and suppose that the 
perturbation $V$ of $H_0$ has the form $V=G^*JG$,
where the operator $G$ is such that $G(\abs{H_0}+I)^{-1/2}$
is compact, and the operator $J=J^*$ is bounded 
and has a bounded inverse
(in contradistinction to \cite{GM},
we do not assume that $J^2=I$; this generalisation is
completely trivial, but may be useful in applications).
Under these assumptions, one can define the perturbed operator 
$H=H_0+G^*JG$.
If $H_0$ is semibounded from below, the sum $H_0+G^*JG$ 
is understood in the form sense.
If $H_0$ is not semibounded from below, one can still define
the operator $H$ using the resolvent identity; 
this is explained in \S\ref{subsection.b.2} below.

Next, we fix an open interval $\d\subset\R$
and assume that the operator $GE_{H_0}(\d)$ belongs to 
the Hilbert--Schmidt class $\SS_2$.
The above assumptions guarantee that for a.e. $\l\in\R$,
the limits $A(\l+i0)$, $B(\l+i0)$ (see \eqref{a.6} or,
for a rigorous definition, \eqref{b.3})
exist in the operator norm and $B(\l+i0)\in\SS_1$.
This implies that the r.h.s. of \eqref{a.7}
(and of its generalisation \eqref{a.10} below)
is well defined.

Further, we accept the following assumption on the 
function $f$ (this assumption will depend on the 
spectral parameter $\l$). 

\begin{condition}
\label{condition.a.1}
Let $\Omega\subset \R$ be a Borel set, and let $f:\Omega\to\R$ 
satisfy the following two conditions at the point $\l$:

(i) 
$\l$ is an interior point of $\Omega$, 
$f$ is continuous and differentiable at $\l$, 
and $f'(\l)>0$;

(ii) 
$\inf\{\abs{f(x)-f(\l)}\mid x\in\Omega,\quad \abs{x-\l}>\d\}>0$
for any $\d>0$.
\end{condition}

We suppose that $\s(H_0)\cup\s(H)\subset\Omega$, the inclusion 
\eqref{a.8} holds and the \cnd{a.1} holds for all $\l\in\d$.
Thus, the SSF for the pair $f(H_0)$, $f(H)$ is well defined. 
Under these assumptions, we prove that for a.e. $\l\in\d$
one has 
\begin{equation}
\xi(f(\l);f(H),f(H_0))=\frac1\pi\int_{-\infty}^\infty\frac{dt}{1+t^2}
\iindex\bigl(E_{\J+A(\l+i0)+tB(\l+i0)}((-\infty,0)),
E_{\J}((-\infty,0))\bigr).
\label{a.10}
%................................................................(a.10)
\end{equation}

In applications to differential operators, 
the hypothesis of the above described result 
(for a suitable choice of the function $f$)
can be easily deduced from the appropriate assumptions 
on the coefficients of the differential operators 
$H_0$, $H$.

{\bf 5.}
Let us describe the idea of the proof. As a main tool, we use a certain
function $\mu(\t;\l,H,H_0)$.
This function is integer valued and depends on two variables 
$\t\in(0,2\pi)$ and $\l\in\R$ and a pair of operators $H_0$, $H$.
The function $\mu$ is closely related to the 
scattering matrix for the pair $H_0$, $H$.
The definition of $\mu$ does not require any trace class assumptions.
However, in the framework of the trace class theory, it is related to 
the SSF for the pair $H_0$, $H$.

In order to define the function $\mu$, we introduce two assumptions 
(\cnd{d.1} and \cnd{d.3}) on the pair $H_0$, $H$. 
The first assumption is formulated in terms of the 
difference of the resolvents of $H_0$ and $H$;
roughly speaking, this assumption means that  $H$ is 
obtained from $H_0$ by means of a relatively compact 
(in an appropriate sense) perturbation.
The other assumption depends on the spectral parameter $\l$ and is 
close to the requirement of the existence of the limits
\eqref{a.6} in the operator norm.

Under these two assumptions, we define the function
$\mu(\t;\l,H,H_0)$ as a {\em spectral flow} 
of a certain family of unitary operators, which depends on 
$H_0$, $H$ and $\l$. The notion of a spectral flow of a 
family of unitary operators is introduced and discussed in 
\S\ref{section.c}. We postpone the definition of 
$\mu$ till \S\ref{section.d}; below we 
only list some of the properties  of $\mu$
(without giving precise statements)
and explain how formula 
\eqref{a.10} can be deduced from these properties.

(i)
Up to an integer constant, $\mu(\t)$ coincides with the eigenvalue
counting function for the spectrum of the scattering matrix
$\calS(\l;H,H_0)$:
\begin{equation}
\mu(\t_1)-\mu(\t_2)=\sum_{\t\in[\t_1,\t_2)}
\dim\Ker(\calS(\l;H,H_0)-e^{i\t}I),
\quad 0<\t_1<\t_2<2\pi.
\label{a.11}
%................................................................(a.11)
\end{equation}

(ii) Suppose that the perturbation $V=H-H_0$ 
can be written down as $V=G^*JG$,
where the operator $G$ is such that $G(\abs{H_0}+I)^{-1/2}$ 
is compact, and $J=J^*$ is bounded and has a bounded inverse.
If the limits \eqref{a.6} exist in the operator norm,
then the following formula for $\mu$ is valid:
\begin{equation}
\mu(\t)=\iindex\bigl(E_{\J}((-\infty,0)),
E_{\J+A(\l+i0)+\cot(\t/2)B(\l+i0)}((-\infty,0))\bigr).
\label{a.14}
%................................................................(a.14)
\end{equation}

(iii) If \eqref{a.1} holds, then $\mu(\t;\l,H,H_0)$ is well defined
for a.e. $\l\in\R$ and the SSF is given by 
\begin{equation}
\xi(\l;H,H_0)=-\frac1{2\pi}\int_0^{2\pi}\mu(\t;\l,H,H_0)d\t.
\label{a.12}
%................................................................(a.12)
\end{equation}

(iv) The function $\mu$ obeys the invariance principle:
\begin{equation}
\mu(\t;\l,H,H_0)=\mu(\t;f(\l),f(H),f(H_0)).
\label{a.13}
%................................................................(a.13)
\end{equation}

Combining \eqref{a.12} and \eqref{a.14} and performing 
the change of variable $t=\cot(\t/2)$ in the resulting integral, 
we obtain \eqref{a.7} (this can be considered as an alternative 
proof of \eqref{a.7}). 
Combining \eqref{a.12}, \eqref{a.13}, \eqref{a.14}, 
we obtain \eqref{a.10}.

Note that, taking into account \eqref{a.11}, the equality \eqref{a.12}
modulo $\Z$ is merely the Birman-Krein formula \eqref{a.5},
and the relation \eqref{a.13} modulo $\Z$ is a trivial consequence 
of the invariance principle for the scattering matrix.
It is the choice of an integer constant that matters in the 
definition of  $\mu$. 
The adequate choice of the constant is related to the normalisation
condition \eqref{a.4}.

In fact, formula \eqref{a.11} is not used in the proof of \eqref{a.10};
we have mentioned it here only in order to explain the 
underlying idea of the proof and the relation between 
the function $\mu$ and the scattering matrix.

{\bf 6. }
Let us describe the structure of the paper.
In  \S\ref{section.b}, we introduce some notation 
and explain in what sense we understand the sum $H_0+G^*JG$
(without assuming that $H_0$ is semibounded from below).
In \S\ref{section.c} we discuss 
the notion of the spectral flow for unitary operators.
In \S\ref{section.d} we define the function $\mu$.
In \S\ref{section.e}, \ref{section.g}, \ref{section.f}, 
we prove formulae \eqref{a.14}, \eqref{a.12}, \eqref{a.13},
respectively.
In \S\ref{section.i}, we prove formula 
\eqref{a.11} and explain the relation of the function 
$\mu$ to the eigenvalue counting functions of the operators
$H_0$, $H$ away from their essential spectrum.

In each section, the  statement and discussion of  
all the results  are given first 
and the proofs are postponed till the end of the section.

{\bf 7. }
In different parts of the paper, we use two different 
points of view on the pair of operators $H_0$, $H$
(in accord with the nature of the question under consideration).
The first point of view is that the `basic' operators 
are the unperturbed operator $H_0$ and the perturbation
$G^*JG$; the perturbed operator $H$ is defined as the sum 
$H=H_0+G^*JG$.
This point of view is aimed at applications. 

According to the second point of view,
the operators $H_0$ and $H$ are defined independently one of 
another and have equal roles; in this case we do not use 
the factorisation of the perturbation $H-H_0$.

%*************************************************************
%*************************************************************
\section{Notation and preliminaries}
\label{section.b}
%************************************************************
%************************************************************

%=============================================================
\subsection{Notation}
\label{subsection.b.1}
%============================================================

{\bf 1. }
Below $\H$, $\K$ are separable Hilbert spaces; 
$I$ is the identity operator.
For a closable linear operator $T:\H\to\K$, 
by $\Dom T$ we denote its domain and 
by $\overline{T}$ --- the closure of $T$.  
For a self-adjoint operator $A$ in a Hilbert space,
the symbols $\s(A)$, $\s_{ess}(A)$, 
$\rho(A)$ denote its spectrum, essential spectrum
and resolvent set
and $E_A(\d)$ is the spectral projection
associated to
a Borel set $\d\subset\R$. 
We also denote by $\Xi(A)$ the $\Xi$ operator associated with $A$
(see \cite{GMN,GM}): $\Xi(A):=E_A((-\infty,0))$.

By ${\mathcal B}(\H,\K)$ we denote the Banach space of all 
bounded operators acting from $\H$ to $\K$;
$\SS_\infty(\H,\K)\subset {\mathcal B}(\H,\K)$ 
is the space of all compact operators and
$\SS_p(\H,\K)$, $p\geq1$, is the standard 
Schatten--von Neumann class.
We write $\calB(\H):=\calB(\H,\H)$,
$\SS_p(\H):=\SS_p(\H,\H)$; 
the norm in the classes $\calB$, $\SS_p$ is denoted by 
$\norm{\cdot}$, 
$\snorm{\cdot}{p}$ and the limits --- by 
$\blim$, $\splim$, respectively.

We shall often use the well-known fact that 
\begin{equation}
A\in\SS_p,\quad 
M_n\xrightarrow{s}0
\quad\Longrightarrow\quad
\|M_nA\|_{\SS_p}\to0,\quad p\in[1,\infty];
\label{d.13}
%..................................................................(d.13)
\end{equation}
here $\xrightarrow{s}$ denotes strong convergence.
If, in addition, $M_n^*\xrightarrow{s}0$, then also 
$\|AM_n\|_{\SS_p}\to0$.
In particular, \eqref{d.13} implies that 
\begin{equation}
A_n\in\SS_p,\quad \|A_n-A\|_{\SS_p}\to0,\quad
M_n\xrightarrow{s}M\quad
\Longrightarrow\quad
\|M_n A_n -MA\|_{\SS_p}\to0.
\label{d.13aa}
%..................................................................(d.13aa)
\end{equation}

Formulas and statements with double indices 
($\pm$ and $\mp$) should be
read as pairs of statements, in one of which 
all the indices take upper
values and in another --- the lower ones.
A constant which first appears
in formula $(i.j)$ is denoted by $C_{i.j}$.
We denote $\C_+=\{z\in\C\mid \Im z>0\}$.
The open ball in a metric space with the centre $x$ 
and radius $r$ is denoted by $B(x;r)$.

{\bf 2.}
A pair $P,Q$ of orthogonal projections in $\H$
is called Fredholm if 
$$
\{+1,-1\}\cap\s_{ess}(P-Q)=\emptyset.
$$
In particular, if $P-Q$ is compact, then the pair 
$P,Q$ is Fredholm.
The index of a Fredholm pair is determined by the formula
\begin{equation}
\iindex(P,Q):=\dim(\Ker(P-Q-I))-
\dim(\Ker(P-Q+I)).
\label{b.1}
%............................................................(b.1)
\end{equation}
Clearly,
$$
\iindex(P,Q)=-\iindex(Q,P).
$$
If either $(P-Q)$ or $(Q-R)$ is 
compact and both $P,Q$ and $Q,R$ are 
Fredholm pairs, then the pair $P,R$ is also Fredholm
and the following chain rule is valid:
\begin{equation}
\iindex(P,R)=\iindex(P,Q)+\iindex(Q,R).
\label{b.1a}
%............................................................(b.1a)
\end{equation}
See, e.g., \cite{AvronSeilerSimon} for the details.

%===================================================================
\subsection{Operator $H(H_0,G,J)$}
\label{subsection.b.2}
%===================================================================

Let $\H$ be a `basic' and $\K$ an `auxiliary' Hilbert space.
Fix a self-adjoint operator  $H_0$ in $\H$ and let 
$G:\H\to\K$ and $J$ in $\K$ be such operators that 
\begin{equation}
\Dom (|H_0|+I)^{1/2}\subset\Dom G,\quad
G(|H_0|+I)^{-1/2}\in\SS_\infty(\H,\K),\quad
J=J^*\in\calB(\K),\quad 0\in\r(J).
\label{b.2}
%..................................................................(b.2)
\end{equation}
Below we define a self-adjoint operator $H$, 
which corresponds to the formal sum
$H_0+G^*JG$. Sometimes we shall explicitly
indicate the dependence of $H$ on $H_0$, $G$, $J$ by
writing $H(H_0,G,J)$.
The construction below goes back to \cite{Kato1}
and is discussed in detail in \cite[\S1.9, 1.10]{Yafaev}.

For $z\in\r(H_0)$ define the following operators 
of the class $\SS_\infty(\K)$:
\begin{equation}
\begin{split}
T(z)&=T(z;H_0,G)=(G(|H_0|+I)^{-1/2})
\frac{|H_0|+I}{H_0-zI}(G(|H_0|+I)^{-1/2})^*,\\
A(z)&=A(z;H_0,G)=\Re T(z),\qquad B(z)=B(z;H_0,G)=\Im T(z).
\end{split}
\label{b.3}
%..................................................................(b.3)
\end{equation}
It is easy to check (see, e.g., \cite[Lemma 1.10.5]{Yafaev}) that
\begin{equation}
0\in\r(I+JT(z))\text{ for all }z\in\C\setminus\R.
\label{b.4}
%..................................................................(b.4)
\end{equation}
Under the assumptions \eqref{b.2}, there exists a unique 
self-adjoint operator  $H=H(H_0,G,J)$ (see \cite[\S1.9, 1.10]{Yafaev}), 
such that for all $z\in\C\setminus\R$ 
its resolvent satisfies the equation
\begin{equation}
(H-zI)^{-1}-(H_0-zI)^{-1}=
-(G(H_0-\overline{z}I)^{-1})^*
(I+JT(z))^{-1}(JG(H_0-zI)^{-1}).
\label{b.5}
%..................................................................(b.5)
\end{equation}
The inverse operator $(I+JT(z))^{-1}$
in the r.h.s. of \eqref{b.5} exists by \eqref{b.4}.
Note that \eqref{b.4} implies 
\begin{equation}
0\in\r(\J+T(z)),\quad z\in\C\setminus\R,
\label{b.7}
%..................................................................(b.7)
\end{equation}
and \eqref{b.5} can be written down as
\begin{equation}
(H-zI)^{-1}-(H_0-zI)^{-1}
=-(G(H_0-\overline{z}I)^{-1})^*
(\J+T(z))^{-1}(G(H_0-zI)^{-1}).
\label{b.8}
%..................................................................(b.8)
\end{equation}

If $H_0$ is semibounded from below, then $H$ coincides with 
the sum $H_0+G^*JG$ in the  form sense.
More precisely, if $h_0[\cdot,\cdot]$ is the sesquilinear form of $H_0$
with the domain $d[h_0](=\Dom (|H_0|+I)^{1/2})$,
then the sesquilinear form $h[\cdot,\cdot]$ of $H$ is defined on the 
domain $d[h]=d[h_0]$ by the relation
$$
h[f,g]=h_0[f,g]+(JGf,Gg),\quad
f,g\in d[h_0].
$$
If the operator $G^*JG$ is well defined and 
$H_0$-bounded with a relative bound $<1$, then 
$H=H_0+G^*JG$ in the sense of the Kato--Rellich theorem.

Finally, by \eqref{b.8}, the difference of the resolvents 
of $H$ and $H_0$ is compact, and therefore the essential spectra 
of $H_0$ and $H$ coincide.

%*****************************************************************
%*****************************************************************
\section{The spectral flow for unitary operators}
\label{section.c}
%****************************************************************
%****************************************************************

%===============================================================
\subsection{Introduction}
\label{subsection.c.1}
%===============================================================
Let $A(t)$, $t\in [0,1]$, be a family of self-adjoint Fredholm operators.
If $A(t)$ is continuous in $t$ in some appropriate sense,
one can define the {\it spectral flow } of $A$, $\spflow(A)$.
A `naive' definition of the spectral flow is the following:
\begin{equation*}
\begin{split} 
\spflow(A)&=\langle\text{the number of eigenvalues of $A(t)$ that
cross $0$ rightwards}\rangle
\\
&\quad-\langle\text{the number of eigenvalues of
$A(t)$ that cross $0$ leftwards}\rangle 
\end{split}
%...........................................................(c.1)
\end{equation*}
as $t$ grows monotonically
from $0$ to $1$. The spectral flow was introduced in 
\cite[\S7]{APS} as the intersection number of the graph
$\cup_{t\in[0,1]}\s(A(t))$ of the spectrum of $A(t)$
with the line $\l=-\e$, where $\e$ is a sufficiently small 
positive number (one can take $\e=0$ if both $A(0)$ and $A(1)$ 
are invertible).
The spectral flow is an important homotopy invariant of the 
family $A(t)$ --- see, e.g., recent treatments in
\cite{RobbinSalamon} and \cite{FPR} and references therein.

In this paper, we will need the notion of the spectral flow 
for {\it unitary}, rather than self-adjoint, operators.
Namely, let us 
fix a Hilbert space $\H$ and a parameter $p\in[1,\infty]$.
Let $Y_p=Y_p(\H)$ be the set of all unitary operators $W$ in $\H$ 
such that $W-I\in \SS_p(\H)$. 
Clearly, $Y_p$ is a metric space with the metric
$d(W_1,W_2)=\snorm{W_1-W_2}{p}$, $p<\infty$ and
$d(W_1,W_2)=\norm{W_1-W_2}$, $p=\infty$.
Consider a mapping $U:[0,1]\to Y_p$. 
We do not suppose that $U$ is continuous; instead, 
we assume that the spectrum $\s(U(t))$ depends continuously 
on $t$ in a certain precise sense to be defined below.
In this section we 
define the spectral flow of the family $U(t)$ 
through the points $z\in\T\setminus\{1\}$.
A `naive' definition of the spectral flow is the following:
\begin{equation}
\begin{split} 
\spflow(z;U)&=\langle\text{the number of eigenvalues of $U(t)$ that
cross $z$ anti-clockwise}\rangle
\\
&\quad-\langle\text{the number of eigenvalues of
$U(t)$ that cross $z$ clockwise}\rangle 
\end{split}
\label{c.1}
%...........................................................(c.1)
\end{equation}
as $t$ grows monotonically
from $0$ to $1$. 

In our subsequent construction, we will have to deal with
$\spflow(z;U)$ as the function of the spectral parameter 
$z\in\T\setminus\{1\}$. For example, we will have to consider the 
integral
$$
\int_0^{2\pi}\spflow(e^{i\theta};U)d\theta
$$
for the families $U:[0,1]\to Y_1$.
Therefore, the behaviour of $\spflow(e^{i\theta};U)$ as an
element of the functional spaces on $(0,2\pi)$ 
(such as $L_1(0,2\pi)$) is essential for us.

Because of this, we find it convenient to give our own 
definition of the spectral flow (see \dfn{c.5} below), 
which is adapted to our 
specific purposes and consistently takes into account 
the dependence of $\spflow(z;U)$ on the spectral parameter
$z$. 

In \S\ref{subsection.c.4} we will show 
that our definition coincides with the 
naive definition  \eqref{c.1} (whenever the latter makes sense)
and therefore is consistent with the standard definition 
of the spectral flow. 
However, in the rest of the paper we do not use this fact 
and work entirely in terms of our definition.

Note that, in contrast to \cite{APS,RobbinSalamon},
our definition does not use the notion of
intersection number and other `difficult' 
topological tools.
We only need the notion of covering space 
(we recall the definition and basic properties of the covering spaces
in \S\ref{subsection.c.1a}).

For the proofs of the main results of this paper we shall 
need only the cases $p=1$, $p=\infty$. 
Nevertheless, we find it instructive to give a universal treatment
of all the cases $p\in[1,\infty]$,
since this does not require any 
considerable modification of the proofs.

%===============================================================
\subsection{Covering spaces}
\label{subsection.c.1a}
%===============================================================
For the reader's convenience, we recall the definition of covering 
spaces and their basic properties. 
The details can be found in any textbook in algebraic topology; 
see, e.g., \cite[Chapter 5]{Massey}.

Let $X$ and $\wt X$ be topological spaces. We suppose that 
$X$ and $\wt X$ are {\it arcwise connected } (i.e., any two
points can be joint by a path)  and {\it locally
arcwise connected } (i.e., any point has a basic family 
of arcwise connected neighbourhoods).
A continuous mapping $\pi:\wt X\to X$ is called a {\it covering, }
if every point $x\in X$ has an arcwise connected open neighbourhood $U$
with the following property. The restriction of $\pi$ onto
each arc component $V$ of $\pi^{-1}(U)$ is a homeomorphism 
between $V$ and $U$.

The important property of covering spaces is that paths and their 
homotopies can be lifted from $X$ to $\wt X$. More precisely:
\begin{proposition}
\label{p.c.0a}
Let $\wt x\in \wt X$, $x=\pi(\wt x)$.
For any path $\g:[0,1]\to X$ with the initial point 
$\g(0)=x$, there exists a unique path (a {\it lift } of $\g$)
 $\wt \g:[0,1]\to \wt X$
such that $\pi\circ\wt\g=\g$ and $\wt\g(0)=\wt x$.
\end{proposition}

The idea of the proof is to express the path $\g$ as a sequence 
of a finite number of `short' paths, each of which is contained 
in an elementary neighbourhood, and then lift each of these paths. 
For the details (and the proof of the uniqueness part), see,
e.g., \cite[Chapter 5, \S3]{Massey}.

\begin{proposition}
\label{p.c.0b}
Let $\wt \g_0,\wt\g_1:[0,1]\to\wt X$ be paths in $\wt X$ 
which have the same initial point: $\wt\g_0(0)=\wt\g_1(0)$.
If $\pi\circ \wt\g_0$ is homotopic to $\pi\circ \wt\g_1$,
then $\wt\g_0$ is homotopic to $\wt\g_1$; in particular, 
$\wt\g_0(1)=\wt\g_1(1)$.
\end{proposition}
The idea of the proof is essentially the 
same as that of \prp{c.0a}.
Let $F:[0,1]\times[0,1]\to X$ be a homotopy between 
$\pi\circ \wt\g_0$ and $\pi\circ \wt\g_1$:
\begin{align*}
F(t,0)&=\pi(\wt\g_0(t)),\quad
&F(t,1)=\pi(\wt\g_1(t)),\\
F(0,s)&=\pi(\wt\g_0(0)),\quad
&F(1,s)=\pi(\wt\g_0(1)).
\end{align*}
Then the square $[0,1]\times[0,1]$ can be subdivided into `small' 
rectangles such that $F$ maps each rectangle into 
an elementary neighbourhood. After that, $F$ can be lifted to 
$\wt X$ locally on each rectangle. The result of this lifting
gives a homotopy between $\wt \g_0$ and $\wt \g_1$.
For the details, see, e.g., \cite[Chapter 5, Lemma 3.3]{Massey}.

%===================================================================
\subsection{The covering $\pi_p:\xtp\to X_p$}
\label{subsection.c.2}
%===================================================================

{\bf 1.} First we define the function space 
$\wt X_p$ which the function $\spflow(\cdot;U)$
will belong to.
Let $\widetilde X_\infty$ be the set of all functions
$f:\T\setminus\{1\}\to\Z$ such that the function
$(0,2\pi)\ni\theta\mapsto f(e^{i\theta})$ is left continuous and
non-increasing.
Clearly, the points $z\in\T\setminus\{1\}$ 
where $f\in\widetilde X_\infty$ is discontinuous, 
can accumulate only to $1$. 
For any $f\in\widetilde X_\infty$, let us introduce the
function $\nu(\cdot;f):\Z\to[0,2\pi]$ by
\begin{equation}
\nu(n;f):=\sup(\{0\}\cup\{\theta\in(0,2\pi)\mid
f(e^{i\theta})>n\}).
\label{c.1a}
%...........................................................(c.1a)
\end{equation}
Clearly, $\nu(\cdot;f)$ is non-increasing
and $\lim_{n\to+\infty}\nu(n;f)=0$,
$\lim_{n\to-\infty}\nu(n;f)=2\pi$.
Note that $f$ can be recovered from $\nu(\cdot;f)$ by the 
formula
\begin{equation}
\label{c.1b}
%...........................................................(c.1b)
f(e^{i\t}):=\inf\{n\in\Z\mid \nu(n;f)<\theta\}.
\end{equation}

For $p\in[1,\infty)$,
let $\wt X_p\subset\wt X_\infty$ be the set of functions
$f$ such that 
$$ 
\sum_{n\geq0}(\nu(n;f))^p+
\sum_{n<0}(2\pi-\nu(n;f))^p<\infty. 
$$ 
For any $p\in[1,\infty]$
and any $f,g\in\xtp$, define 
$$
\rtp(f,g):=\|\nu(\cdot;f)-\nu(\cdot;g)\|_{l_p(\Z)}. 
$$
Note that 
$$
\wt\r_1(f,g)=\int_0^{2\pi}\lvert f(e^{i\t})-g(e^{i\t})\rvert d\t.
$$

\begin{proposition}
\label{p.c.1}
The function $\rtp$ is a metric on $\xtp$.
With respect to this metric, $\wt X_p$ is arcwise connected
and locally arcwise connected.
\end{proposition}

{\bf 2.} Consider the following equivalence relation on $\xtp$: 
$$ 
f\sim g\Longleftrightarrow \exists n\in\Z: 
\forall z\in\T\setminus \{1\}, 
\quad f(z)=g(z)+n. 
$$ 
Let $X_p$ be the quotient space $\xtp/\!\!\!\sim$, and let
$\pi_p:\xtp\to X_p$ be the corresponding projection.
For $f,g\in X_p$ define 
$$ 
\rho_p(f,g)=\inf\{\rtp(\wt f,\wt g)\mid
\pi_p(\wt f)=f, \pi_p(\wt g)=g\}. 
$$

\begin{proposition}
\label{p.c.2}
The function $\rho_p$ is a metric on $X_p$.
With respect to this metric, $X_p$ is arcwise connected
and locally arcwise connected.
\end{proposition}

Obviously, the mapping $\pi_p:\xtp\to X_p$ is continuous.

\begin{proposition}
\label{p.c.3}
The mapping $\pi_p:\xtp\to X_p$ is a covering.
\end{proposition}

Clearly, an element $f\in X_p$ is uniquely determined by
specifying the set of discontinuities $z_n\in\T\setminus\{1\}$
of an element $\wt f\in\pi_p^{-1}(f)$ 
together with the heights $m(z_n)$ of the `jumps' of
$\wt f$ at the points $z_n$.
Thus,  the space $X_p$ can be identified with 
the set of the spectra of 
all unitary operators $W\in Y_p$;
under this identification, $z_n$ become eigenvalues
with the multiplicities $m(z_n)$.

{\em Notation}
Let $\g:[0,1]\to \wt X_p$ be any mapping.
Then $\g$ depends on two variables, $t\in[0,1]$ 
and $z\in\T\setminus\{1\}$. If we need to indicate the 
dependence of $\g$ on both variables $z$ and $t$,
we write $\g(z;t)$.
If $\g$ is considered as an element of the function space $\wt X_p$
(for a fixed $t$), we write $\g(t)$.

{\bf 3.}
It is obvious that the following diagram is commutative for any
$1\leq q<r\leq\infty$:
\begin{equation}
\begin{CD}
\wt X_q@>in_{\wt X_q\to\wt X_r}>> \wt X_r\\
@V\pi_q VV @VV\pi_r V\\
X_q@>in_{X_q\to X_r}>>X_r
\end{CD}
\label{c.4}
%............................................................(c.4)
\end{equation}
Here $in_{\wt X_q\to\wt X_r}$
and $in_{X_q\to X_r}$ are the natural embeddings.

%=================================================================
\subsection{The mapping $\eta_p:Y_p\to X_p$}
\label{subsection.c.3}
%=================================================================

{\bf 1.} 
Below we use the following natural notation for the arcs 
of the unit circle on the complex plane:
$$
(e^{i\t_1},e^{i\t_2})=\{e^{i\t}\mid\t_1<\t<\t_2\},
\quad\t_1<\t_2,
$$
with the obvious modifications for
$[e^{i\t_1},e^{i\t_2}]$,
$(e^{i\t_1},e^{i\t_2}]$,
$[e^{i\t_1},e^{i\t_2})$.

Let $W\in Y_p$ and $\theta_1,\theta_2\in(0,2\pi)$. Define
\begin{equation}
N(e^{i\theta_1}, e^{i\theta_2};W)=\left\{
\begin{array}{lc}
\rank E_W([e^{i\theta_1},e^{i\theta_2})),&\theta_1<\theta_2,\\
0,&\theta_1=\theta_2,\\ -\rank E_W
([e^{i\theta_2},e^{i\theta_1})),&\theta_2<\theta_1.
\end{array}
\right.
\label{c.4a}
%............................................................(c.4a)
\end{equation}

It is easy to see that for any $z_0\in\T\setminus\{1\}$
the function $\T\setminus\{1\}\ni z\mapsto N(z,z_0;W)\in\Z$ 
belongs to the space $\xtp$.
\begin{proposition}
\label{p.c.3a}
Fix $z_0\in\T\setminus\{1\}$. The mapping 
$$ 
Y_p\ni W\mapsto
N(\cdot,z_0;W)\in \xtp
$$ 
is continuous at the `points' $W$ such that 
$z_0\in\T\setminus\s(W)$.
\end{proposition}

{\bf 2.} 
Let us define the mapping $\eta_p$:
\begin{equation}
Y_p\ni W\mapsto \eta_p(W):=\pi_p(N(\cdot,z_0;W))\in X_p,\quad
z_0\in\T\setminus\s(W).
\label{c.6}
%............................................................(c.6)
\end{equation}
Clearly, this definition does not depend on $z_0$, since the
change of $z_0$ results in adding an integer constant to
$N(\cdot,z_0;W)$. By \prp{c.3a}, the mapping $\eta_p$ 
is continuous.

{\bf 3.}
Note that the following diagram is commutative for any
$1\leq q<r\leq\infty$:
\begin{equation}
\begin{CD}
Y_q@>in_{Y_q\to Y_r}>> Y_r\\
@V\eta_q VV @VV\eta_r V\\
X_q@>in_{X_q\to X_r}>>X_r
\end{CD}
\label{c.7}
%............................................................(c.7)
\end{equation}
Here $in_{X_q\to X_r}$ and $in_{Y_q\to Y_r}$ 
are the natural embeddings.

%=================================================================
\subsection{The spectral flow}
\label{subsection.c.4}
%=================================================================

{\bf 1.}
Now we are ready to define the spectral flow of a family 
$U:[0,1]\to Y_p$.
But first we have to take into account 
one complication of a formal nature.
In our construction below (see \S\ref{d.1}) we have to deal with 
the families, defined on an open, rather than closed, interval
$(0,1)$. 
At the same time, it appears that the composition
$\eta_p\circ U$ can be extended by continuity to the endpoints
$\pm1$. 
Thus, first we need the notation for such an extension. 
Suppose that a mapping $\g:(0,1)\to X_p$ is continuous and
the limits $\lim_{t\to0+}\g(t)$, $\lim_{t\to1-}\g(t)$ exist. 
Then we write that {\em the extension of $\g$ exists} and denote by  
$$
\ext(\g)
$$ 
the mapping $\g$, extended by continuity to the whole 
interval $[0,1]$. 
\begin{definition}
\label{d.c.5}
Let $U:(0,1)\to Y_p$ be such a  mapping that 
the extension $\g:=\ext(\eta_p\circ U)$ exists.
Let $\wt \g$ be a lift of  $\g$ into $\wt X_p$.
Then we define 
\begin{equation}
\spflow(z;U):=\wt\g(z;1)-\wt\g(z;0).
\label{c.7a}
%............................................................(c.7a)
\end{equation}
\end{definition}

{\it  \dfn{c.5}
does not depend on the choice 
of the lift $\wt\g$. } 
Indeed, let $\wt\g_1$ and $\wt\g_2$ be 
two lifts of $\g$.
Then the function
$\wt\g_2(0)-\wt\g_1(0)$ is an integer constant;
let us denote this constant by $n$.
By the uniqueness of the lift of a path 
with a fixed initial point, one has
$\wt\g_2(t)\equiv\wt\g_1(t)+n$ and therefore
$\wt\g_2(1)-\wt\g_2(0)=
\wt\g_1(1)-\wt\g_1(0)$.

{\it \dfn{c.5} does not depend on $p$ 
in the following sense. }
Let $1\leq q<r\leq\infty$ and 
let $U_q:(0,1)\to Y_q$ be such a mapping that the extension
$\g_q=\ext(\eta_q\circ U_q)$ exists.
Let $\wt\g_q$ be the lift of $\g_q$ and 
$\wt\g_q(1)-\wt\g_q(0)$ be the spectral flow of $U_q$.

Further, consider the mapping 
$U_r:=in_{Y_q\to Y_r}\circ U_q:(0,1)\to Y_r$.
It follows from  \eqref{c.7} that the extension
$\g_r=\ext(\eta_r\circ U_r)$ exists and 
$\g_r=in_{X_q\to X_r}\circ \g_q$.
Consider the lift $\wt\g_r$ of $\g_r$.
Taking into account \eqref{c.4},
one sees that $in_{\wt X_q\to\wt X_r}\circ\wt\g_q$
is also a lift of $\g_r$. 
From here it follows that 
$$
in_{\wt X_q\to\wt X_r}(\wt\g_q(1))-
in_{\wt X_q\to\wt X_r}(\wt\g_q(0))=
\wt\g_r(1)-\wt\g_r(0).
$$

{\bf 2.} Thus defined, the spectral flow is homotopy invariant:
\begin{proposition}
\label{p.c.5}
Let $U_1, U_2:(0,1)\to Y_p$ be two mappings 
such that the extensions  $\g_1=\ext(\eta_p\circ U_1)$
and $\g_2=\ext(\eta_p\circ U_2)$ exist and are homotopic
(in particular, this implies that $\g_1(0)=\g_2(0)$ and 
$\g_1(1)=\g_2(1)$).
Then 
\begin{equation}
\spflow(z;U_1)=\spflow(z;U_2),\quad z\in\T\setminus\{1\}.
\label{c.7c}
%............................................................(c.7c)
\end{equation}
\end{proposition}
\begin{proof} A direct application of \prp{c.0b}.
\end{proof}

Note that our proof of the invariance principle
\eqref{a.13} depends heavily on the homotopy invariance of the 
spectral flow.

{\bf 3.}
In this paper we do not explicitly use the fact that 
\dfn{c.5} agrees with the `naive'
definition \eqref{c.1}, whenever the latter makes sense.
However, let us give a sketch of proof of this fact.
Here for the sake of  simplicity of notation we assume that our 
mappings $U$ are already defined on the whole of $[0,1]$
and thus need not be extended.

First suppose that for a mapping 
$U:[0,1]\to Y_p$ (such that $\eta_p\circ U$ is continuous)
there exists $z_0\in\T\setminus\{1\}$ such that 
$z_0\in\r(U(t))$ for all $t\in[0,1]$.
One easily checks that in this case,
according to \dfn{c.5},
$$
\spflow(z;U)=N(z,z_0;U(1))-N(z,z_0;U(0)).
$$
Clearly, this agrees with \eqref{c.1}.

Further, for an arbitrary mapping $U:[0,1]\to Y_p$ 
(such that $\eta_p\circ U$ is continuous),
one can always find a finite cover of $[0,1]$ 
by the intervals $\d_n$, $n=1,\dots,N$, 
with the property that for any $n$ there exists
$z_n\in\T\setminus\{1\}$, $z_n\in\r(U(t))$ 
for any $t\in\d_n$.
In this case, one can write 
\begin{equation}
\spflow(z;U)=\sum_{n=1}^N
\bigl(N(z,z_n;U(t_n))-N(z,z_n;U(t_{n-1}))\bigr)
\label{c.7b}
%............................................................(c.7b)
\end{equation}
for a set of points $0=t_0<t_1<\dots<t_N=1$,
$t_n\in\d_n\cap\d_{n+1}$ for $n=1,\dots,N-1$.
Formula \eqref{c.7b} also agrees with \eqref{c.1}.

%================================================================
\subsection{Proof of Propositions \ref{p.c.1}---\ref{p.c.3a}}
\label{subsection.c.5}
%================================================================

{\bf 1.} {\it Proof of \prp{c.1}} 
1. Let us prove that $\wt \r_p$ is a metric.
Clearly, $\rtp(f,g)=\rtp(g,f)$ and
$\rtp(f,g)\geq0$. Suppose that $f\not\equiv g$;
by \eqref{c.1b}, it follows that 
$\nu(\cdot;f)\not\equiv\nu(\cdot;g)$
and therefore $\wt\r_p(f,g)\not=0$.

The triangle inequality for $\rtp$ is evident. 

2. We shall prove that any ball in $\wt X_p$ is arcwise connected;
clearly, this will imply that $\wt X_p$ is arcwise connected and 
locally  arcwise connected. 

For every $f_0,f_1\in\wt X_p$, let
$$
\nu_\a(n)=\a\nu(n;f_1)+(1-\a)\nu(n;f_0),
\quad\a\in[0,1],\quad n\in\Z.
$$
The formula \eqref{c.1b}
recovers the family $f_\a$ of the functions such that 
$\nu(n;f_\a)=\nu_\a(n)$.
Clearly,  the path $[0,1]\ni\a\mapsto f_\a\in\wt X_p$ 
connects $f_0$ and $f_1$; moreover,
$\wt\r_p(f_0,f_\a)\leq\wt\r_p(f_0,f_1)$.
Thus, every ball in $\wt X_p$ is arcwise connected.
\qed

{\bf 2.} {\it Auxiliary facts} 

1.
Note that 
\begin{equation}
\rtp(f+n,g+n)=\rtp(f,g) 
\mbox{ for any constant }n\in\Z. 
\label{c.3}
%............................................................(c.3)
\end{equation}

2. 
Clearly, for any $f\in\xtp$ one has
\begin{equation}
\inf_{n\in\Z\setminus\{0\}}\rtp(f+n,f)=\rtp(f+1,f)>0. 
\label{c.9}
%............................................................(c.9)
\end{equation}

3.
Let us prove that
\begin{equation}
\forall f,g\in\xtp\quad \exists n\in\Z:\quad
\inf_{m\in\Z}\rtp(f+m,g)=\rtp(f+n,g). 
\label{c.10}
%............................................................(c.10)
\end{equation}
In other words, the infimum in \eqref{c.10} is always attained.

First let $p\not=\infty$. Then, clearly, 
$$
\lim_{|m|\to\infty}\wt\r_p(f+m,g)=\infty, 
$$ 
which proves \eqref{c.10}.
Next, let $p=\infty$. Then 
$$
\lim_{|m|\to\infty}\wt\rho_\infty(f+m,g)=2\pi, 
$$ 
whereas $\wt\rho_\infty(f+m,g)\leq2\pi$ for any $m$. This proves
\eqref{c.10} for $p=\infty$.

{\bf 3.} {\it Proof of \prp{c.2}}
1. Let us prove that $\r_p$ is a metric.
Clearly, $\r_p(f,g)=\r_p(g,f)$ and
$\r_p(f,g)\geq0$. Suppose that $\r_p(f,g)=0$; let us check that
$f=g$. Fix $\wt f\in\pi_p^{-1}(f)$, $\wt g\in\pi_p^{-1}(g)$. By
\eqref{c.10}, the relation $\r_p(f,g)=0$ implies that $\rtp(\wt
f+n,\wt g)=0$ for some $n\in\Z$ and thus $\wt f+ n=\wt g$ and
therefore $f=g$.

The triangle inequality for $\r_p$ follows directly from the 
triangle inequality for $\wt\r_p$.

2. Obviously, $\pi_p(\wt X_p)=X_p$. Since $\wt X_p$ is arcwise connected,
it follows that $X_p$ is also arcwise connected.

3. Let us prove that $X_p$ is locally arcwise connected.
To this end, we prove that every ball in $X_p$ is arcwise connected.
Fix $f\in X_p$, $\wt f\in\pi_p^{-1}(f)$ and $r>0$ and consider 
the open ball $B(f;r)$ with the centre $f$ and radius $r$.
Below we prove that $\pi_p$ maps the ball $B(\wt f;r)$ onto
$B(f;r)$. Since $B(\wt f;r)$ is arcwise connected (see 
the proof of \prp{c.1}), this will imply that $B(f;r)$ is 
also arcwise connected.

The inclusion $\pi_p(B(\wt f;r))\subset B(f;r)$ is evident.
Let us prove that $B(f;r)\subset\pi_p(B(\wt f;r))$.
If $g\in B(f;r)$ and $\wt g\in\pi_p^{-1}(g)$,
then $\inf_{m\in\Z}\wt\r_p(\wt f+m,\wt g)<r$, 
which, by \eqref{c.10},
implies that $\wt\r_p(\wt f+m,\wt g)<r$
for some $m\in\Z$.
Thus, $\wt\r_p(\wt f,\wt g-m)<r$ and therefore 
$\wt g-m\in B(\wt f;r)$ and $g=\pi_p(\wt g-m)\in\pi_p(B(\wt f;r))$.
\qed

{\bf 4.} {\it Proof of \prp{c.3} } 
Fix $f\in X_p$, $\wt f\in\pi_p^{-1}(f)$ and $\e<\rtp(\wt f+1,\wt f)/3$. 
Let us prove that the ball $B(f;\e)$ is an elementary neighbourhood.
We shall prove that 
$\pi_p^{-1}(B(f;\e))=\cup_{n\in\Z}B(\wt f+n;\e)$, where
the balls $B(\wt f+n;\e)$ are mutually disjoint, arcwise connected and 
the restriction $\pi_p\mid B(\wt f+n;\e)$ is a homeomorphism between 
$B(\wt f+n;\e)$ and $B(f;\e)$.

Let us first check that the balls  $B(\wt f+n;\e)$ are mutually disjoint. 
Indeed, let $\wt g\in B(\wt f+n;\e)\cap B(\wt f+m;\e)$. 
Then 
$\rtp(\wt f+n,\wt f+m)\leq\rtp(\wt f+n,\wt g)+\rtp(\wt g,\wt f+m)<2\e$.
By \eqref{c.9} and the choice of $\e$, the last inequality implies $m=n$.

In the course of the proof of \prp{c.2}, we have checked that 
$\pi_p(B(\wt f+n;\e))=B(f;\e)$ for any $n\in\Z$.
The same reasoning also shows that 
$\pi_p^{-1}(B(f;\e))=\cup_{n\in\Z}B(\wt f+n;\e)$.

Let us prove that the restriction
$\pi_p\mid B(\wt f+n;\e)$ is injective.
Let $\pi_p (\wt g)=\pi_p(\wt h)$ for $\wt g, \wt h\in B(\wt f+n;\e)$.
Then $\wt g=\wt h+m$ for some $m\in\Z$. 
Using \eqref{c.3}, one has:
\begin{equation*}
\begin{split}
\wt g\in  B(\wt f+n;\e)
\quad&\Rightarrow\quad\rtp(\wt f+n,\wt g)<\e\quad 
\Rightarrow\quad
\rtp(\wt f+n-m,\wt h)<\e\quad\Rightarrow\quad\wt h\in 
 B(\wt f+n-m;\e)
\\
&\Rightarrow\quad m=0\quad\Rightarrow\quad \wt g=\wt h. 
\end{split}
\end{equation*}

4. Finally, let us check that $(\pi_p\mid B(\wt f+n;\e))^{-1}$ is continuous.
Let $\wt g, \wt h\in B(\wt f+n;\e)$, $g=\pi_p(\wt g)$, $h=\pi_p(\wt h)$.
Below we show that if $\r_p(g,h)<\e$, then $\rtp(\wt g,\wt h)=\r_p(g,h)$.
Indeed, by \eqref{c.10}, one has $\r_p(g,h)=\rtp(\wt g+m,\wt h)$ 
for some $m\in\Z$.
Let us show that $m=0$. Using \eqref{c.3}, one has
\begin{equation*}
\begin{split}
\rtp(\wt f+m,\wt f)&=\rtp(\wt f+n+m,\wt f+n)\leq\rtp(\wt f+n+m,\wt g+m)+
\rtp(\wt g+m,\wt h)
\\
&\quad+\rtp(\wt h,\wt f+n)<3\e,
\end{split}
\end{equation*}
which, by \eqref{c.9} and the
choice of $\e$, implies $m=0$. \qed

{\bf 5.} The proof of \prp{c.3a} is based on the following 
\begin{lemma}
\label{l.c.4}
For any $\e\in(0,2\pi)$ there exists $C_{\ref{c.5}}(\e)>0$ 
such that for any $z_0\in\T\setminus\{1\}$ 
and any operators $W_1$, $W_2\in Y_p$ with the property 
$$
[z_0e^{-i\e},z_0e^{i\e}]\cap\s(W_j)=\emptyset,
\quad j=1,2,
$$
the following estimate holds:
\begin{equation}
\rtp(N(\cdot,z_0;W_1),N(\cdot,z_0;W_2))
\leq C_{\ref{c.5}}(\e)\snorm{W_1-W_2}{p}.
\label{c.5}
%............................................................(c.5)
\end{equation}
\end{lemma}

1. Let us first prove the following auxiliary statement. 
For an operator $A=A^*\in\SS_p$, 
let $\{\l^{(+)}_n(A)\}_{n\in\N}$ be the sequence of its 
non-negative eigenvalues 
listed in decreasing order counting multiplicities, and let
$\l^{(-)}_n(A):=\l^{(+)}_n(-A)$.
Denote $\Zo=\Z\setminus\{0\}$. 
Let $\L(A)\in l_p(\Zo)$ be the sequence
$$
\L_n(A)=\left\{
\begin{array}{ll}
\l_n^{(+)}(A),& n>0;\\
\l_{-n}^{(-)}(A),& n<0.
\end{array}
\right.
$$
Let us prove that for any self-adjoint operators $A_1,A_2\in\SS_p$,
\begin{equation}
\|\L(A_1)-\L(A_2)\|_{l_p(\Zo)}\leq\|A_1-A_2\|_{\SS_p}.
%..........................................................{c.11}
\label{c.11}
\end{equation}
For $p=\infty$, the above relation follows directly from the 
variational characterisation of the eigenvalues. 
The general case is a consequence of 
a slight modification of Lidski's theorem 
\cite{Lidski} (see also \cite[Chapter 2, \S6.5]{Kato}).
First note that it is sufficient to prove \eqref{c.11} for
finite rank operators $A_1$, $A_2$. 
In the finite rank case, Lidski's theorem says that 
\begin{equation}
\l_n(A_1)-\l_n(A_2)=\sum_m\s_{nm}\l_m(A_1-A_2),
%..........................................................{c.12}
\label{c.12}
\end{equation}
where $\{\l_n(A)\}$ is the sequence of all (positive and negative)
eigenvalues of $A$, listed in the order of decreasing of 
the absolute value $|\l_n(A)|$, 
and $\s_{nm}$ is a matrix satisfying
\begin{equation}
\sum_n|\s_{nm}|\leq1,\quad
\sum_m|\s_{nm}|\leq1.
%..........................................................{c.13}
\label{c.13}
\end{equation}
The relations \eqref{c.12}, \eqref{c.13} imply 
(cf. \cite{Kato}) that
$$
\sum_n|\l_n(A_1)-\l_n(A_2)|^p\leq\sum_n|\l_n(A_1-A_2)|^p=
\|A_1-A_2\|^p_{\SS_p},\quad p\in[1,\infty),
$$
which differs from the desired inequality \eqref{c.11}
only by the method of numbering the eigenvalues. 
Following the proof of Lidski's theorem, it is 
not difficult to see that it holds also in the case 
when the positive and negative eigenvalues are numbered 
separately; more precisely, one has
\begin{equation}
\begin{split}
\l_n^{(\pm)}(A_1)-\l_n^{(\pm)}(A_2)&=
\sum_m\s_{nm}^{(\pm)}\l_m(A_1-A_2),\\
\sum_n\abs{\s_{nm}^{(+)}}+\abs{\s_{nm}^{(-)}}
&\leq1,\quad
\sum_m\abs{\s_{nm}^{(\pm)}}\leq1.
\end{split}
%..........................................................{c.14}
\label{c.14}
\end{equation}
In the same way as above, \eqref{c.14} implies \eqref{c.11}.

2. Below we will need the following fact. 
For any $\f\in C^\infty(\T)$ and any two 
unitary operators $W_1$, $W_2$ such that 
$W_1-W_2\in\SS_p$, one has 
\begin{equation}
\|\f(W_1)-\f(W_2)\|_{\SS_p}\leq C_{\ref{c.18}}(\f)\|W_1-W_2\|_{\SS_p}.
\label{c.18}
%..........................................................{c.18}
\end{equation}
In order to prove \eqref{c.18}
(see, e.g., \cite[\S5.4]{BYa} for the details 
and discussion), one first writes a representation
$$
\f(z)=\sum_{n\in\Z}c_nz^n,\quad \sum_{n\in\Z}|n||c_n|<\infty,
$$
which is valid for all smooth enough $\f$. 
Next, it is easy to check that 
$$
\|W_1^n-W_2^n\|_{\SS_p}\leq n\|W_1-W_2\|_{\SS_p}.
$$
Therefore, \eqref{c.18} holds with 
$C_{\ref{c.18}}(\f)=\sum_{n\in\Z}|n||c_n|$.

3. Now we are ready to prove the estimate  \eqref{c.5}.
Let $\f_\e\in C^\infty(\T)$ be such a function that
$\f_\e(e^{i\theta})=\theta$ for all
$\theta\in[-2\pi+\e,-\e]$.
Denote $\f_{\e,z_0}(z):=\f_\e(z/z_0)+\arg z_0$,
where $\arg z_0\in(0,2\pi)$.
It is straightforward to see that for $j=1,2$
and $n=1,2,\dots$, one has
\begin{gather}
\nu(n-1;N(\cdot,z_0;W_j))=\l^{(+)}_{n}(\f_{\e,z_0}(W_j)),\quad
\nu(-n;N(\cdot,z_0;W_j))=2\pi-\l^{(-)}_{n}(\f_{\e,z_0}(W_j)),
\notag\\
\intertext{and therefore}
\rtp(N(\cdot,z_0;W_1),N(\cdot,z_0;W_2))=
\|\L(\f_{\e,z_0}(W_1))-\L(\f_{\e,z_0}(W_2))\|_{l_p(\Zo)}.
\label{c.16}
%..........................................................{c.16}
\end{gather}
The relations \eqref{c.16}, \eqref{c.11} and \eqref{c.18} 
together imply \eqref{c.5} with the constant 
$$
C_{\ref{c.5}}(\e)=
\sup_{z_0\in\T\setminus\{1\}}C_{\ref{c.18}}(\f_{\e,z_0}).
\quad \qed
$$

{\em Proof of \prp{c.3a}}
Fix $W_0$ such that $z_0\in\T\setminus\s(W_0)$
and $\e>0$ such that 
$[z_0e^{-i\e},z_0e^{i\e}]\cap\s(W_0)=\emptyset$.
Then for any $W\in Y_p$ such that $\norm{W-W_0}<\e/2$,
one has 
$[z_0e^{-i\e/2},z_0e^{i\e/2}]\cap\s(W)=\emptyset$.
Thus, we can apply \lma{c.4}, which yields
$$
\rtp(N(\cdot,z_0;W),N(\cdot,z_0;W_0))\leq
C_{\ref{c.5}}(\e/2)\snorm{W-W_0}{p}.
$$
Clearly, this implies the continuity of the mapping in hand at the 
`point' $W_0$.
\qed

%================================================================
\subsection{Lemma on convergence in $X_p$}
\label{subsection.c.6}
%================================================================
In the proof of \thm{f.2} below we shall need the following
\begin{lemma}
\label{l.c.5}
Let $W_n$ and $W'_n$ be sequences of operators in
$Y_p$ such that $\limn\snorm{W_n-W'_n}{p}=0$.
Then the limit 
$\xplim_{n\to\infty}\eta_p(W_n)$ exists if and only if
the limit $\xplim_{n\to\infty}\eta_p(W'_n)$ exists.
If these limits exist, they coincide.
\end{lemma}
{\em Proof}
1. For any $f\in X_\infty$, let us introduce the notation 
$$
\s(f):=\{\exp(i\nu(n;\wt f))\mid n\in\Z\}\cup\{1\},\quad
\wt f\in\pi_\infty^{-1}(f)
$$
(recall that $\nu(n;\wt f)$ is defined by \eqref{c.1a}).
Clearly, this definition does not depend on the choice of an
element $\wt f\in\pi_\infty^{-1}(f)$.
It is also clear that in this notation,
$$
\s(W)=\s(\eta_\infty(W)),\quad W\in Y_\infty.
$$

2. Suppose that the limit 
$f:=\xplim_{n\to\infty}\eta_p(W_n)$ exists.
Below we prove that the limit 
$\xplim_{n\to\infty}\eta_p(W'_n)$ also exists 
and is equal to $f$.
Fix $z_0\in\T\setminus\s(f)$ and $\e>0$ such that 
$[z_0e^{-i\e},z_0e^{i\e}]\cap\s(f)=\emptyset$.
If $n$ is large enough so that 
$\r_\infty(f,\eta_\infty(W_n))<\e/3$, we get 
$$
[z_0e^{-i2\e/3},z_0e^{i2\e/3}]\cap\s(\eta_\infty(W_n))=\emptyset.
$$
Further, if $n$ is large enough so that 
$\r_\infty(f,\eta_\infty(W_n))<\e/3$
and $\norm{W_n-W'_n}<\e/3$, we get 
$$
[z_0e^{-i\e/3},z_0e^{i\e/3}]\cap\s(\eta_\infty(W'_n))=\emptyset.
$$
For such $n$ we can apply \lma{c.4}, which yields
$$
\r_p(\eta_p(W_n),\eta_p(W'_n))\leq
C_{\ref{c.5}}(\e/3)\snorm{W_n-W'_n}{p}\to0
\quad\text{as $n\to\infty$.}
$$
Thus, 
$\limn\r_p(\eta_p(W'_n),f)=0$.
\qed

%***************************************************************
\section{The function $\mu$: definition}
\label{section.d}
%***************************************************************

%=============================================================
\subsection{Definition}
\label{subsection.d.1}
%=============================================================

Let $H_0$ and $H$ be self-adjoint operators in a Hilbert space $\H$.
For any $z\in\r(H_0)\cap\r(H)$ define a unitary operator in $\H$ by 
\begin{equation}
M(z;H,H_0):=
\frac{H-\overline{z}I}{H-zI}\frac{H_0-zI}{H_0-\overline{z}I}=
(I+(z-\overline{z})(H-zI)^{-1})
(I+(\overline{z}-z)(H_0-\overline{z}I)^{-1}).
\label{d.1}
%..................................................................(d.1)
\end{equation}
Next, in what follows we fix  $p\in[1,\infty]$. We introduce 
\begin{condition}
\label{condition.d.1}
(i)
For any $z\in\r(H_0)\cap\r(H)$ one has 
\begin{equation}
(H-zI)^{-1}-(H_0-zI)^{-1}\in\SS_p.
\label{d.8aa}
%...............................................................(d.8aa)
\end{equation}
(ii)
For any $\l\in\R$ one has
\begin{equation}
\lim_{y\to+\infty}y
\snorm{(H-(\l+iy)I)^{-1}-(H_0-(\l+iy)I)^{-1}}{p}=0.
\label{d.3a}
%..................................................................(d.3a)
\end{equation}
\end{condition}
By the identity 
\begin{equation}
M(z)-I=(z-\overline{z})
((H-zI)^{-1}-(H_0-zI)^{-1})
\frac{H_0-zI}{H_0-\overline{z}I},
\label{d.4}
%..................................................................(d.4)
\end{equation}
the inclusion \eqref{d.8aa} is equivalent to 
\begin{equation}
M(z;H,H_0)-I\in \SS_p(\H),
\label{d.3}
%..................................................................(d.3)
\end{equation}
and the relation \eqref{d.3a} is equivalent to 
\begin{equation}
\lim_{y\to+\infty}\snorm{M(\l+iy;H,H_0)-I}{p}=0.
\label{d.4a}
%..................................................................(d.4a)
\end{equation}
\begin{proposition}
\label{p.d.1a}
(i) If \eqref{d.8aa} holds for one value of $z$,
then it holds for all $z\in\r(H_0)\cap\r(H)$.

(ii)
If \eqref{d.3a} holds for one value of $\l$, then 
it holds for all $\l\in\R$.

(iii) \cnd{d.1}(i) implies that the mapping
$$
\C\setminus\R\ni z\mapsto M(z;H,H_0)-I\in \SS_p(\H)
$$
is continuous. 
\end{proposition}
Further, we need one more assumption.
Recall that the class $X_p$ and the mapping 
$\eta_p$ have been defined 
in \S\ref{subsection.c.2}, \ref{subsection.c.3}.
Fix $\l\in\R$.
\begin{condition}
\label{condition.d.3}
The limit 
\begin{equation}
\xplim_{y\to0+}\eta_p(M(\l+iy;H,H_0))
\label{d.6}
%..................................................................(d.6)
\end{equation}
exists.
\end{condition}
Under the Assumptions \ref{condition.d.1} and \ref{condition.d.3},
consider the mapping 
\begin{equation}
U:(0,1)\ni t\mapsto M(\l+i(1-t)t^{-1};H,H_0)\in Y_p.
\label{d.7}
%..................................................................(d.7)
\end{equation}
Clearly, the mapping $U$ satisfies the hypothesis of \dfn{c.5}
and therefore $\spflow(z;U)$ is well defined.

\begin{definition}
\label{definition.d.4}
Suppose that for a pair of selfadjoint operators $H_0$, $H$
and for $\l\in\R$, 
the Assumptions \ref{condition.d.1}, \ref{condition.d.3}
hold true.
Let $U$ be the mapping \eqref{d.7}; then we define
\begin{equation}
\mu(\theta;\l,H,H_0):=\spflow(e^{i\t};U),
\quad\theta\in(0,2\pi).
\label{d.8}
%..................................................................(d.8)
\end{equation}
\end{definition}

%======================================================================
\subsection{Sufficient conditions}
\label{subsection.d.2}
%======================================================================

Let $\H$ be a `basic' and $\K$ an `auxiliary' Hilbert spaces and let 
operators $H_0$, $G$, $J$, $H=H(H_0,G,J)$ be as described in
\S\ref{subsection.b.2}.
Below we give sufficient conditions (in terms of $H_0$, $G$, $J$), 
which ensure that the Assumptions \ref{condition.d.1} and 
\ref{condition.d.3} hold true for the pair $H_0$, $H$.
In addition to \eqref{b.2}, assume that 
\begin{equation}
G(|H_0|+I)^{-1/2}\in\SS_{2p}(\H,\K)
\label{d.8a}
%...............................................................(d.8a)
\end{equation}
for some $p\in[1,\infty]$.
\begin{proposition}
\label{p.d.5}
Assume \eqref{b.2}, \eqref{d.8a}.
Then, for the pair of operators $H_0$, $H$, 
\cnd{d.1} holds true.
\end{proposition}

\begin{proposition}
\label{p.d.6}
Assume \eqref{b.2}, \eqref{d.8a}
and define the operators \eqref{b.3}.
Suppose that for some $\l\in\R$ 

{\rm (i)} 
the limit $\slim_{y\to0+}(\J +T(\l+iy))^{-1}$ exists;

{\rm (ii)}
the limit 
$\splim_{y\to0+}B(\l+iy)=:B(\l+i0)$ exists.

Then, for the pair $H_0$, $H$, 
\cnd{d.3} holds at the point $\l$.
\end{proposition}

\begin{proposition}
\label{p.d.7}
Assume \eqref{b.2}, \eqref{d.8a} 
and suppose that for an open 
interval $\d\subset\R$ one has 
\begin{equation}
GE_{H_0}(\d)\in\SS_2(\H,\K).
\label{d.9}
%..................................................................(d.9)
\end{equation}
Then for a.e. $\l\in\d$

{\rm (i)} 
the  limits
\begin{equation}
\sqlim_{y\to0+}T(\l+iy), \quad
\splim_{y\to0+}B(\l+iy)
\label{d.9a}
%...................................................................(d.9a)
\end{equation}
exist, where $q=p$ if $p>1$ and $q$ is any number greater than $1$,
if $p=1$;

{\rm (ii)} one has $0\in\r(\J+T(\l+i0))$.

Thus, the hypotheses (i), (ii) of 
\prp{d.6} hold true and the pair $H_0$, $H$ 
satisfies \cnd{d.3}.
\end{proposition}

%=====================================================================
\subsection{Operator $S(z)$}
\label{subsection.d.3}
%=====================================================================
In order to prove Propositions \ref{p.d.5}--\ref{p.d.7},
below we introduce an auxiliary operator $S(z)$. 
Let $\H$ be a `basic' and $\K$ an `auxiliary' Hilbert spaces.
Let the operators $H_0$, $G$, $J$
be as described in \S\ref{subsection.b.2};
assume \eqref{b.2} and \eqref{d.8a} for some $p\in[1,\infty]$
and let $H=H(H_0,G,J)$.
For any $z\in\C\setminus\R$ define 
\begin{equation}
S(z)=S(z;H_0,G,J):=
I-2i B^{1/2}(z)(\J+T(z))^{-1}B^{1/2}(z).
\label{d.10}
%..................................................................(d.10)
\end{equation}
The inverse operator in the r.h.s. of \eqref{d.10} exists by \eqref{b.7}.
A straightforward calculation shows that $S(z)$ is unitary in $\K$.
Clearly, $S(z)-I\in \SS_p$. The operator $S(z)$ can also be presented as 
\begin{equation*}
\begin{split}
S(z)&= I-2iB^{1/2}(z)(I+JT(z))^{-1}JB^{1/2}(z)=\\
&=I-2iB^{1/2}(z)J(I+T(z)J)^{-1}B^{1/2}(z).
\end{split}
\end{equation*}
The definition of the operator $S(z)$ 
copies the stationary representation for the scattering matrix 
(see \eqref{i.1}). For this reason, the operators
of this type are well studied (see, e.g., 
\cite{BYa2} and references therein).

\begin{lemma}
\label{l.d.8}
Assume \eqref{b.2} and \eqref{d.8a}. 
Then the mapping 
\begin{equation}
\C\setminus\R\ni z\mapsto S(z)-I\in \SS_p(\K)
%..................................................................(d.11)
\label{d.11}
\end{equation}
is continuous and 
\begin{equation}
\|S(z)-I\|_{\SS_p}\to0\text{ as }\Im z\to+\infty.
\label{d.12}
%..................................................................(d.12)
\end{equation}
\end{lemma}
{\it Proof}

1. Let us first check that  
\begin{equation}
\text{the mapping}\quad
\r(H_0)\ni z\mapsto T(z)\in \SS_p \quad\text{is continuous}
\label{d.13a}
%..................................................................(d.13a)
\end{equation}
and 
\begin{equation}
\|T(z)\|_{\SS_p}\to0\quad\text{ as }\Im z\to+\infty.
\label{d.13b}
%..................................................................(d.13b)
\end{equation}
In order to do this, observe that the mapping 
\begin{equation}
\r(H_0)\ni z\mapsto \frac{|H_0|+I}{H_0-zI}
\in{\mathcal B}(\H)
\label{d.14}
%..................................................................(d.14)
\end{equation}
is continuous (in the operator norm) and
\begin{equation}
\frac{|H_0|+I}{H_0-zI}\xrightarrow{s}0\text{ as }
\Im z\to+\infty.
\label{d.14a}
%..................................................................(d.14a)
\end{equation}

Now recall the definition \eqref{b.3} of $T(z)$.
By \eqref{d.13}, the relation \eqref{d.13a} follows from  
\eqref{d.8a} and the continuity of \eqref{d.14}.
Similarly, \eqref{d.13b} follows from 
\eqref{d.8a} and \eqref{d.14a}.

2. Clearly, the relations \eqref{d.13a} and \eqref{b.7} 
imply that
\begin{equation}
\text{the mapping }\quad\C\setminus\R\ni z
\mapsto (\J+T(z))^{-1}\in\calB(\H)\quad\text{ is continuous.}
\label{d.14b}
%..................................................................(d.14b)
\end{equation}

3. By \eqref{d.13aa}, the relations \eqref{d.13a} and 
\eqref{d.14b}
imply the continuity of the mapping \eqref{d.11}. 
The relation \eqref{d.13b} implies \eqref{d.12}.
\qed

\begin{theorem}
\label{t.d.10}
Assume \eqref{b.2} and let $H=H(H_0,G,J)$. For any $z\in\C\setminus\R$
the operator \linebreak $M(z;H,H_0)-I$ is compact and 
\begin{equation}
\eta_\infty(M(z;H,H_0))=\eta_\infty(S(z;H_0,G,J)).
\label{d.16}
%..................................................................(d.16)
\end{equation}
\end{theorem}
{\it Proof} 
1. 
By \eqref{d.4} and \eqref{b.8}, one has 
\begin{multline}
M(z)=I-(z-\overline{z})(G(H_0-\overline{z}I)^{-1})^*\\
\times(\J+T(z))^{-1}(G(H_0-zI)^{-1})
(I-(z-\overline{z})(H_0-\overline{z}I)^{-1}).
\label{d.17}
%..................................................................(d.17)
\end{multline}
It follows that $M(z)-I\in\SS_\infty$.

2. For $R>0$, denote $P^{(R)}=E_{H_0}((-R,R))$, 
$G^{(R)}=GP^{(R)}$, $H_0^{(R)}=H_0P^{(R)}$.
Note that $G^{(R)}\in\SS_\infty(\H,\K)$ and 
$H_0^{(R)}\in{\mathcal B}(\H)$.
Further, let 
$H^{(R)}=H_0^{(R)}+(G^{(R)})^*JG^{(R)}(\in{\mathcal B}(\H))$.
By \eqref{d.13}, the relation 
$P^{(R)}=(P^{(R)})^*\xrightarrow{s}I$ implies that 
$$
\|G^{(R)}(|H_0|+I)^{-1/2}-G(|H_0|+I)^{-1/2}\|\to0\text{ as }
R\to+\infty,
$$
and thus 
$$
\|T(z;H_0^{(R)},G^{(R)})-T(z;H_0,G)\|\to0\text{ as }R\to+\infty.
$$
By the definition \eqref{d.10} of $S(z)$ it follows that 
$$
\|S(z;H_0^{(R)},G^{(R)},J)-S(z;H_0,G,J)\|\to0\text{ as }R\to+\infty
$$
and by \eqref{d.17} it follows that 
$$
\|M(z;H^{(R)},H_0^{(R)})-M(z;H,H_0)\|\to0
\text{ as }R\to+\infty.
$$
Therefore, since the mapping 
$\eta_\infty:Y_\infty\to X_\infty$ is continuous,
it is sufficient to prove that 
\begin{equation}
\eta_\infty(M(z;H^{(R)},H_0^{(R)})=\eta_\infty(S(z;H_0^{(R)},G^{(R)},J))
\label{d.18}
%..................................................................(d.18)
\end{equation}
for any $R>0$. For the sake of brevity, below we suppress the index
$R$ in the notation and suppose that $H_0\in{\mathcal B}(\H)$ and
$G\in\SS_\infty(\H,\K)$. 
We also denote $V:=G^*JG$.

3. Recall that for any two bounded operators $A$, $B$ 
and any $\l\not=0$ one has 
\begin{equation}
\dim\Ker(AB-\l I)=\dim\Ker(BA-\l I).
\label{d.18a}
%..................................................................(d.18a)
\end{equation}

By \eqref{d.18a}, for any $\l\not=1$, one has
\begin{equation*}
\begin{split}
\dim\Ker(M(z)-\l I)&=
\dim\Ker\br{\frac{H-\overline{z}I}{H-zI}\frac{H_0-zI}{H_0-\overline{z}I}
-\l I}\\
&=\dim\Ker\br{(H-\overline{z}I)(H_0-\overline{z}I)^{-1}
((H-zI)(H_0-zI)^{-1})^{-1}-\l I}
\\
&=\dim\Ker\br{(I+V(H_0-\overline{z}I)^{-1})(I+V(H_0-zI)^{-1})^{-1}
-\l I}\\
&=\dim\Ker\br{I-2iV\Im((H_0-zI)^{-1})(I+V(H_0-zI)^{-1})^{-1}
-\l I}\\
&=\dim\Ker\br{I-2iG\Im((H_0-zI)^{-1})(I+V(H_0-zI)^{-1})^{-1}G^*J
-\l I}.
\end{split}
\end{equation*}
A direct computation shows that 
$$
(I+V(H_0-zI)^{-1})^{-1}G^*J=
G^*(\J+T(z))^{-1}.
$$
Thus,
\begin{equation*}
\begin{split}
\dim\Ker\br{M(z)-\l I}&= 
\dim\Ker\br{I-2iG\Im ((H_0-zI)^{-1})G^*(\J+T(z))^{-1}-\l I}\\
&=\dim\Ker\br{I-2iB(z)(\J+T(z))^{-1}-\l I}\\
&=\dim\Ker\br{I-2iB^{1/2}(z)(\J+T(z))^{-1}B^{1/2}(z)-\l I}\\
&=\dim\Ker\br{S(z)-\l I},
\end{split}
\end{equation*}
which implies \eqref{d.16}.
\qed

%=================================================================
\subsection{Proofs of Propositions \protect\ref{p.d.1a},
\protect\ref{p.d.5}--\protect\ref{p.d.7}}
\label{subsection.d.4}
%=================================================================

{\it Proof of \prp{d.1a} } 
(i) follows from the identity
\begin{equation}
(H-zI)^{-1}-(H_0-zI)^{-1}=
\frac{H-z_0I}{H-zI}((H-z_0I)^{-1}-(H_0-z_0I)^{-1})
\frac{H_0-z_0I}{H_0-zI}.
%..................................................................(d.18b)
\label{d.18b}
\end{equation}

(ii) Suppose that \eqref{d.3a} holds for $\l=\l_0$.
In \eqref{d.18b}, take $z=\l+iy$, $z_0=\l_0+iy$.
Now the desired assertion follows from the fact that 
$$
\sup_{y>1}\norm{\frac{H-(\l_0+iy)I}{H-(\l+iy)I}}<\infty,
\quad
\sup_{y>1}\norm{\frac{H_0-(\l_0+iy)I}{H_0-(\l+iy)I}}<\infty.
$$

(iii)
Let us use \eqref{d.17} and check that the r.h.s. of 
this identity depends continuously on $z$ in the $\SS_p$ norm.
Similarly to the proof of \lma{d.8}, factorising
$$
G(H_0-zI)^{-1}=[G(|H_0|+I)^{-1/2}]
[(|H_0|+I)^{1/2}(H_0-zI)^{-1}],
$$
and using \eqref{d.13}, we check that 
the operator $G(H_0-zI)^{-1}$ depends continuously 
on $z$ in $\SS_{2p}$ norm.
Taking into account \eqref{d.14b} and the fact that the operator 
$(I-(z-\overline{z})(H_0-\overline{z}I)^{-1})$
depends continuously on $z$ in the operator norm, 
we get the desired assertion.
\qed

{\it Proof of \prp{d.5}}
Let us use \eqref{b.8}. Since $(\J+T(z))^{-1}$ is bounded and 
$G(H_0-zI)^{-1}\in\SS_{2p}$, we get the inclusion \eqref{d.8aa}.
The relation \eqref{d.3a} is equivalent to \eqref{d.4a};
the latter follows from \thm{d.10} and \eqref{d.12}.
\qed

{\it Proof of \prp{d.6} }
By \thm{d.10}, it is sufficient to prove that the limit 
$$
\splim_{y\to0+}(S(\l+iy;H_0,G,J)-I)
$$
exists. By \eqref{d.13aa}, the existence of the above limit 
follows directly from the definition of operator $S$ 
and the hypothesis of the proposition.
\qed

{\it Proof of \prp{d.7}}

1. For any $\d'\subset\R$, denote 
$$
T_{\d'}(z)=T(z;H_0,GE_{H_0}(\d')).
$$
Denoting $\D=\R\setminus\d$, we see that 
$$
T(z)=T_\d(z)+T_\D(z).
$$
It is one of the classical results of 
the trace class scattering theory (see \cite{BE,Naboko})
that the inclusion \eqref{d.9} implies that 
for a.e. $\l\in\R$ the limit $T_\d(\l+i0)$ exists 
in $\SS_r(\K)$ (for any $r>1$) and the limit 
$\lim_{y \to0+}\Im T_\d(\l+iy)$ exists in $\SS_1(\K)$. 
On the other hand, the function $T_\D(z)\in\SS_p(\K)$ is analytic
in $\C\setminus\D$ and $\Im T_\D(\l)=0$ for all $\l\in\d$. 
Thus, for a.e. $\l\in\d$ the limits
\eqref{d.9a} exist.

2. It remains to check that the limit 
$\blim_{y\to0+}(\J+T(\l+iy))^{-1}$ exists for a.e. $\l\in\d$. 
In order to do this, write 
$$
(\J+T(z))^{-1}=(\J+T_\D(z))^{-1}(I+F(z))^{-1},\quad
F(z)=T_\d(z)(\J+T_\D(z))^{-1}.
$$
Let us check that for a.e. $\l\in\R$ the limits 
\begin{equation}
\blim_{y\to0+}(\J+T_\D(\l+iy))^{-1}\quad\text{ and }\quad
\blim_{y\to0+}(I+F(\l+iy))^{-1}
\label{d.19}
%..................................................................(d.19)
\end{equation}
exist.

3. By the Fredholm analytic alternative, the set 
$$
{\mathcal N}=\{\l\in\d\mid 0\in\s(\J+T_\D(\l))\}
$$
is discrete in $\d$ (i.e., the points of ${\mathcal N}$ 
can possibly accumulate only to the endpoints 
of the interval $\d$).
Thus, the limit 
$\blim_{y\to0+}(\J+T_\D(\l+iy))^{-1}$ exists for all 
$\l\in\d\setminus{\mathcal N}$.

4. The function $F(z)\in\SS_1(\K)$ is analytic in $\C_+$ 
and for a.e. $\l\in\d$ has limit values $F(\l+i0)$ in $\SS_q(\K)$
(for any $q>1$). Thus, using Theorem 1.8.5 from \cite{Yafaev},
we obtain that the limit  $\blim_{y\to0+}(I+F(\l+iy))^{-1}$ 
exists for a.e. $\l\in\d$.
\qed

%**********************************************************
\section{Formula for $\mu$}
\label{section.e}
%*********************************************************

%=======================================================
\subsection{Statement of the result}
\label{subsection.e.1}
%=======================================================
Let the operators $H_0$, $G$, $J$ 
be as described in \S\ref{subsection.b.2}, 
assume \eqref{b.2} and let $H=H(H_0,G,J)$.
Recall that for a self-adjoint operator $A$, we denote 
$\Xi(A):=E_A((-\infty,0))$.
\begin{theorem}
\label{t.e.1}
Suppose that, for some $\l\in\R$,
the limit $\blim_{y\to0+}T(\l+i\e)$ exists 
and $0\in\r(\J+T(\l+i0))$.
Then for all $\theta\in(0,2\pi)$ the pair of projections
$\Xi(\J)$, $\Xi(\J+A(\l+i0)+\cot(\t/2)B(\l+i0))$
is Fredholm and
\begin{equation}
\mu(\theta;\l,H,H_0)=
\iindex\bigl(\Xi(\J), \Xi(\J+A(\l+i0)+\cot(\t/2)B(\l+i0))\bigr).
\label{e.1}
%...............................................................(e.1)
\end{equation}
\end{theorem}

If $J=\pm I$, then \eqref{e.1} takes the form 
\begin{align*}
\mu(\theta;\l,H,H_0)&=
-\rank E_{A(\l+i0)+\cot(\t/2)B(\l+i0)}((-\infty,-1)),
\quad &J=I,\\
\mu(\theta;\l,H,H_0)&=
\rank E_{A(\l+i0)+\cot(\t/2)B(\l+i0)}([1,\infty)),
\quad &J=-I.
\end{align*}
Note that, in particular, this implies the following monotonicity 
rule for the function $\mu$:
$$
\pm J\geq0\quad\Rightarrow\quad
\mp\mu(\t;\l,H,H_0)\geq0.
$$
Related statements are well known in the spectral analysis 
of the scattering matrix --- 
see \cite{BYa2} and references therein.

The relation \eqref{e.1} also implies 
the following estimates for $\mu$:
$$
\pm\mu(\t;\l,H,H_0)\leq\rank \Xi(\pm J).
$$
In particular, if the perturbation $G^*JG$ has rank 
$n<\infty$, then the absolute value of $\mu$ 
does not exceed $n$.

%=======================================================
\subsection{The spectrum of $S(z)$}
\label{subsection.e.2}
%=======================================================
Consider the following operators $A$, $B$, $J$:
\begin{equation}
\begin{split}
A&=A^*\in\SS_{\infty}(\K),\quad
0\leq B\in\SS_\infty(\K),\quad
J=J^*\in{\mathcal B}(\K),\\
0&\in\r(J),\quad
0\in\r(\J+A+iB).
\end{split}
\label{e.2}
%...............................................................(e.2)
\end{equation}
Under these assumptions, define a unitary operator in $\K$ by 
\begin{equation}
S=I-2iB^{1/2}(\J+A+iB)^{-1}B^{1/2}.
\label{e.3}
%...............................................................(e.3)
\end{equation}
The proof of \eqref{e.1} is based on the following simple 
characterisation of the spectrum of $S$.
\begin{lemma}
\label{l.e.2}
Assume \eqref{e.2} and let $S$ be defined by \eqref{e.3}.
Then for any $\theta\in(0,2\pi)$ one has
\begin{equation}
\dim\Ker(S-e^{i\t}I)=\dim\Ker(\J+A+\cot(\t/2)B).
\label{e.4}
%...............................................................(e.4)
\end{equation}
\end{lemma}
{\it Proof } One has (using \eqref{d.18a}):
\begin{equation*}
\begin{split}
\dim\Ker(S-e^{i\t}I)
&=\dim\Ker(I-2iB(\J+A+iB)^{-1}-e^{i\t}I)\\
&=\dim\Ker((\J+A-iB)(\J+A+iB)^{-1}-e^{i\t}I)\\
&=\dim\Ker(\J+A-iB-e^{i\theta}(\J+A+iB))\\
&=\dim\Ker(\J+A+\cot(\theta/2)B).\quad
\qed
\end{split}
\end{equation*}

We shall need the following auxiliary statement, which is 
a very slight modification of one of the results of \cite{GM}.

\begin{lemma}
\label{l.e.3}
Let $M=M^*\in{\mathcal B}(\K)$,
$0\leq B\in\SS_\infty(\K)$ and 
$0\in\r(M+\tau B)$ for some $\tau\in\R$.
Then $\Xi(M),\Xi(M+B)$ is a Fredholm pair of projections and 
\begin{equation}
\iindex(\Xi(M),\Xi(M+B))=
\sum_{s\in(0,1]}\dim\Ker(M+sB).
\label{e.5}
%...............................................................(e.5)
\end{equation}
\end{lemma}

{\it Proof }
1. 
In \cite[Corollary 4.8]{GM}, the desired assertion has been 
proven under the additional assumption $B\in\SS_1(\K)$.
Below we show that this assumption can be lifted.

2. 
First note that the condition $0\in\r(M+\tau B)$ implies that 
$0\not\in\s_{ess}(M)$.
Further, it is easy to see that 
$$
\Xi(M)-\Xi(M+B)\in\SS_\infty(\K).
$$
This can be proven by representing the above projections 
by Riesz integrals and using the resolvent identity
(cf. \cite[Lemmas 3.5, 3.8]{GM}). The above inclusion
implies that $\Xi(M)$, $\Xi(M+B)$ is a Fredholm pair.

3.
First assume that $0\in\r(M)$ and $0\in\r(M+B)$. 
Let $0\leq B_n\in\SS_1(\K)$, $\|B_n-B\|\to0$ as $n\to\infty$.
For all large enough $n$ we will have $0\in\r(M+\tau B_n)$.
By  \cite[Corollary 4.8]{GM}, for such $n$ one has
\begin{equation}
\iindex(\Xi(M),\Xi(M+B_n))=
\sum_{s\in(0,1]}\dim\Ker(M+sB_n).
\label{e.6}
%...............................................................(e.6)
\end{equation}
Our aim is to pass to the limit in \eqref{e.6}.

4. By \cite[Theorem 3.12]{GM}, the l.h.s. of \eqref{e.6} 
tends to the l.h.s. of \eqref{e.5} as $n\to\infty$.
Further, by the Birman--Schwinger principle in a gap 
(see, e.g., \cite{Birman}), one has 
$$
\sum_{s\in(0,1]}\dim\Ker(M+sB)=
\rank E_{B^{1/2}M^{-1}B^{1/2}}((-\infty,-1]).
$$
Since $\|B_n^{1/2}M^{-1}B_n^{1/2}-
B^{1/2}M^{-1}B^{1/2}\|\to0$, we see that the r.h.s of 
\eqref{e.6} tends to the r.h.s. of \eqref{e.5}.

5. 
In order to get rid of the assumptions
$0\in\r(M)$, $0\in\r(M+B)$, we observe that for all 
small enough $\e>0$ one has
$0\in\r(M+\e B)$, $0\in\r(M+B+\e B)$ and thus 
$$
\iindex(\Xi(M+\e B),\Xi(M+B+\e B))=
\sum_{s\in(\e,1+\e]}\dim\Ker(M+sB).
$$
Taking $\e\to0+$ in the above formula, we get \eqref{e.5}.
\qed

\begin{lemma}
\label{l.e.4}
Assume \eqref{e.2} and let $S$ be defined by \eqref{e.3}.
Then for the function $N(\cdot,\cdot;S)$, defined by 
\eqref{c.4a}, one has 
for any $\t_1,\t_2\in(0,2\pi)$:
\begin{equation}
\begin{split}
N(e^{i\t_1},e^{i\t_2};S)&=
\iindex\bigl(\Xi(\J+A+\cot(\t_2/2)B),\Xi(\J+A+\cot(\t_1/2)B)\bigr)\\
&=\iindex\bigl(\Xi(\J),\Xi(\J+A+\cot(\t_1/2)B)\bigr)\\
&\quad+\iindex\bigl(\Xi(\J+A+\cot(\t_2/2)B),\Xi(\J)\bigr);
\end{split}
\label{e.7}
%...............................................................(e.7)
\end{equation}
all the three pairs of projections in the r.h.s. are Fredholm.
\end{lemma}
{\it Proof }
1. 
First of all we note that 
\begin{equation}
\Xi(\J+A+\cot(\t_j/2)B)-\Xi(\J)\in\SS_\infty(\K),
\quad j=1,2.
\label{e.7a}
%.....................................................(e.7a)
\end{equation}
As in the previous lemma, 
this can be proven by representing 
$\Xi(\J+A+\cot(\t_j/2)B)$ and $\Xi(\J)$ 
by the Riesz integrals and using the resolvent identity
(cf. \cite[Lemmas 3.5, 3.8]{GM}).
The inclusion \eqref{e.7a} implies that all the three pairs
of projections in the r.h.s. of \eqref{e.7} are Fredholm.

2.
It is sufficient to prove \eqref{e.7} for $\t_1<\t_2$.
Indeed, the case $\t_1>\t_2$ follows from the above mentioned one 
by changing the roles of $\t_1$ and $\t_2$; 
for $\t_1=\t_2$ the relation \eqref{e.7} trivially holds.

In the case  $\t_1<\t_2$, 
using Lemmas \ref{l.e.2} and \ref{l.e.3}, one has:
\begin{equation*}
\begin{split}
\rank E_S([e^{i\t_1},e^{i\t_2}))&=
\sum_{\t\in[\t_1,\t_2)}\dim\Ker(S-e^{i\t}I)
=\sum_{\t\in[\t_1,\t_2)}\dim\Ker(\J+A+\cot(\t/2)B)\\
&=\sum_{\t\in(\cot(\t_2/2),\cot(\t_1/2)]}\dim\Ker(\J+A+tB)\\
&=\iindex(\Xi(\J+A+\cot(\t_2/2)B),\Xi(\J+A+\cot(\t_1/2)B)).
\end{split}
\end{equation*}
Note  that \lma{e.3} is applicable, since,
by the analytic Fredholm alternative, the assumption
$0\in\r(\J+A+iB)$ (see \eqref{e.2}) implies that 
$0\in\r(\J+A+\tau B)$ for all $\tau\in\R$ but for 
a discrete set of points.

3. 
Thus, we have proven the first equality in \eqref{e.7}. 
The second one follows by the chain rule \eqref{b.1a}. 
Note that  the inclusion \eqref{e.7a} 
ensures the applicability of the chain rule.
\qed

%==========================================================
\subsection{Proof of \protect\thm{e.1}}
\label{subsection.e.3}
%==========================================================
{\bf 1.} First we need a simple result which shows that 
the r.h.s. of \eqref{e.1} depends continuously on
$A(\l+i0)$ and $B(\l+i0)$.
This statement is closely related to
\cite[Lemma 2.5]{Push} and \cite[Theorem 3.12]{GM}.

\begin{lemma}
\label{l.e.5}
Assume \eqref{e.2} and let, in addition,
$B\in\SS_p$, $p\in[1,\infty]$.
Let $A_k=A_k^*\in\SS_\infty(\K)$,
$0\leq B_k\in\SS_p(\K)$, $J_k=J_k^*\in\calB(\K)$,
$k\in\N$ be such operators that $0\in\r(J_k)$, 
$0\in\r(\J_k+A_k+iB_k)$,
$\lim_{k\to\infty}\norm{A_j-A}=0$,  
$\lim_{k\to\infty}\snorm{B_k-B}{p}=0$,
$\lim_{k\to\infty}\norm{J_j-J}=0$.
Define the functions
\begin{align*}
f:\T\setminus\{1\}\ni e^{i\t}&\mapsto
f(e^{i\t})=\iindex(\Xi(\J), \Xi(\J+A+\cot(\t/2)B))\in\Z,\\
f_k:\T\setminus\{1\}\ni e^{i\t}&\mapsto
f_k(e^{i\t})=\iindex(\Xi(\J_k), \Xi(\J_k+A_k+\cot(\t/2)B_k))\in\Z.
\end{align*}
Then $f,f_k\in\wt X_p$ and 
\begin{equation}
\wt \r_p(f_k,f)\to0\quad\text{ as }k\to\infty.
%...............................................................(e.8)
\label{e.8}
\end{equation}
\end{lemma}
{\it Proof }
1. Define the operator $S$ by \eqref{e.3} and let 
$$
S_k=I-2iB_k^{1/2}(\J_k+A_k+iB_k)^{-1}B_k^{1/2}.
$$
As in \prp{d.6}, we see that  
$\snorm{S_k-S}{p}\to0$ as $k\to\infty$.

Fix $\t_0\in(0,2\pi)$ such that $e^{i\t_0}\in\r(S)$.
By \prp{c.3a}, 
\begin{equation}
\wt \r_p(N(\cdot,e^{i\t_0};S_k),N(\cdot,e^{i\t_0};S))\to0
\quad\text{ as }k\to\infty.
%...............................................................(e.9)
\label{e.9}
\end{equation}

2. By \lma{e.4}, 
\begin{align*}
N(e^{i\t},e^{i\t_0};S)&=f(e^{i\t})+C(\t_0),\\
N(e^{i\t},e^{i\t_0};S_k)&=f_k(e^{i\t})+C_k(\t_0)
\end{align*}
with
\begin{align*}
C(\t_0)&=\iindex(\Xi(\J+A+\cot(\t_0/2)B),\Xi(\J)),\\
C_k(\t_0)&=\iindex(\Xi(\J_k+A_k+\cot(\t_0/2)B_k),\Xi(\J_k)).
\end{align*}
Since $e^{i\t_0}\in\r(S)$, by \lma{e.2} 
one has $0\in\r(\J+A+\cot(\t_0/2)B)$.
By \cite[Theorem 3.12]{GM}, it follows that  
$\lim_{k\to\infty}C_k(\t_0)=C(\t_0)$.
Since $C_k(\t_0)$ and $C(\t_0)$ are integer valued, 
one has $C_k(\t_0)=C(\t_0)$ for all large enough $k$.
Thus, by \eqref{c.3}, the relation 
\eqref{e.9} implies \eqref{e.8}.
\qed

{\bf 2.} {\it Proof of \thm{e.1} }
1. First of all, we note that for all $\t\in(0,2\pi)$
$$
\Xi(\J)-\Xi(\J+A(\l+i0)+
\cot(\theta/2)B(\l+i0))\in\SS_\infty(\K)
$$
(cf. \eqref{e.7a}) and thus the 
pair of projections in the r.h.s. of \eqref{e.1} is Fredholm.

2. Let $U$ be the mapping \eqref{d.7} (for $p=\infty$)
and $\g=\ext(\eta_\infty\circ U)$
(remind that $\eta_\infty$ has been introduced  
in \S\ref{subsection.c.3}, and 
$\ext$ --- in \S\ref{subsection.c.4}).
Below we explicitly construct the lift of $\g$.
Let us define the mapping $\wt \g:[0,1]\to\wt X_\infty$ by
\begin{align*}
\wt \g(e^{i\t};0)&=0;\\
\wt\g(e^{i\t};t)&=
\iindex\bigl(\Xi(\J), \Xi(\J+A(z)+\cot(\t/2)B(z))\bigr),
\quad z=\l+i(1-t)t^{-1},
\quad t\in(0,1);\\
\wt\g(e^{i\t};1)&=
\iindex\bigl(\Xi(\J), \Xi(\J+A(\l+i0)+\cot(\t/2)B(\l+i0))\bigr).
\end{align*}
Below we show that:

(i) $\wt \g$ is continuous;

(ii) $\pi_\infty\circ\wt\g=\g$.

The statements (i), (ii)  mean that $\wt\g$ is
the lift of $\g$ with  
$\wt\g(0)=0$. Since the r.h.s. of
\eqref{e.1} coincides with $\wt\g(e^{i\t};1)$,
this implies the statement of the theorem.

3. 
By \lma{e.5}, the continuity of $\wt\g$ 
for $t\in(0,1)$
follows from the norm continuity of $A(z)$, $B(z)$
(see \eqref{d.13a}) in $z$.
Similarly, the continuity of $\wt\g$ at $t=0$ 
follows from \eqref{d.13b} and 
the continuity at $t=1$  is evident.

The relation $\pi_\infty\circ\wt\g=\g$
follows from \thm{d.10} and \lma{e.4}. 
\qed

%*****************************************************************
\section{The function $\mu$ and the perturbation determinant}
\label{section.g}
%*****************************************************************

%===========================================================
\subsection{Statement of the result}
\label{subsection.g.1}
%===========================================================
Let the operators $H_0$, $G$, $J$
be as described in \S\ref{subsection.b.2}.
Assume \eqref{b.2} and  \eqref{d.8a} with $p=1$
and let $H=H(H_0,G,J)$.
As in \cite[\S8.1.4]{Yafaev}, we introduce the 
`modified perturbation determinant' 
\begin{equation}
D_{H/H_0}(z)=\det(I+JT(z)),\quad z\in\r(H_0).
\label{g.1}
%..............................................................(g.1)
\end{equation}
If the operator $V=G^*JG$ is well defined and 
$V(H_0-zI)^{-1}\in\SS_1(\H)$, then $D_{H/H_0}(z)$ 
coincides with the usual perturbation determinant
$\D_{H/H_0}(z)$.
By \eqref{d.13a}, the determinant $D_{H/H_0}(z)$ 
is continuous in $z\in\r(H_0)$ 
(it is, of course, even analytic in $z$, 
but we do not use this fact).
By \eqref{d.13b} with $p=1$, one has $D_{H/H_0}(z)\to0$ 
as $\Im z\to+\infty$.
Let us fix the branch of $\arg D_{H/H_0}(z)$ by 
\begin{equation}
\arg D_{H/H_0}(z)\to0\text{ as }\Im z\to+\infty.
\label{g.2}
%..............................................................(g.2)
\end{equation}

By Propositions \ref{p.d.5}, \ref{p.d.7}, 
for $p=1$ and a.e. $\l\in\R$, the 
Assumptions \ref{condition.d.1}
and \ref{condition.d.3}
hold true.
Therefore, for a.e. $\l\in\R$ the function 
$\mu(\cdot;\l,H,H_0)$ is well defined
and belongs to $L_1(0,2\pi)$.

\begin{theorem}
\label{t.g.1}
Assume \eqref{b.2} and  \eqref{d.8a} with $p=1$,
define the function $D_{H/H_0}$ by \eqref{g.1} 
and fix the branch of $\arg D_{H/H_0}$ by \eqref{g.2}. 
Then for a.e. $\l\in\R$ the limit
$\lim_{y\to0+}\arg D_{H/H_0}(\l+iy)$  exists and
\begin{equation}
\begin{split}
\lim_{y\to0+}\arg D_{H/H_0}(\l+iy)&=
-\frac12\int_0^{2\pi}\mu(\t;\l,H,H_0)d\t\\
&=\int_{-\infty}^{\infty}\frac{dt}{1+t^2}
\iindex\br{\Xi(\J+A(\l+i0)+tB(\l+i0)),\Xi(\J)}.
\end{split}
\label{g.3}
%..............................................................(g.3)
\end{equation}
\end{theorem}

\begin{remark}
A similar reasoning shows that 
under the hypothesis of \thm{g.1}
$$
\arg D_{H/H_0}(z)
=\int_{-\infty}^{\infty}\frac{dt}{1+t^2}
\iindex\br{\Xi(\J+A(z)+tB(z)),\Xi(\J)},
\quad z\in\C_+.
$$
This formula might be of an independent interest, 
although we do not need it in this paper.
\end{remark}

Recalling the Krein's formula \eqref{a.3}, \eqref{a.4} 
for the SSF, we see that for $G\in\SS_2(\H,\K)$,
the first equation in \eqref{g.3} implies 
\eqref{a.12}.
The second equation (in the case $J^2=I$) 
leads us to the representation \eqref{a.7}, 
which was originally obtained in \cite{GM}.

%============================================================
\subsection{Proof of \protect\thm{g.1}}
%============================================================

1. First let us prove that 
\begin{equation}
\det M(z;H,H_0)=\overline{D_{H/H_0}(z)}/D_{H/H_0}(z),
\quad z\in\C_+.
\label{g.5}
%..............................................................(g.5)
\end{equation}
One has:
\begin{align*}
\overline{D_{H/H_0}(z)}/D_{H/H_0}(z)&=
\det\br{(I+JT(\overline{z}))(I+JT(z))^{-1}}\\
&=\det\br{(I+JT(z)-2iJB(z))(I+JT(z))^{-1}}\\
&=\det\br{I-2iJB(z)(I+JT(z))^{-1}}\\
&=\det\br{I-2iB^{1/2}(z)(I+JT(z))^{-1}JB^{1/2}(z)}
=\det S(z;H_0,G,J).
\end{align*}
Finally, note that, by \thm{d.10},
$$
\det S(z;H_0,G,J)=\det M(z;H,H_0).
$$

2. It follows from \eqref{g.5} that 
$$
\arg D_{H/H_0}(z)=-\frac12\arg\det M(z;H,H_0),
$$
where the branches are fixed by \eqref{g.2} 
and by the condition
\begin{equation}
\arg\det M(z;H,H_0)\to 0\text{ as }\Im z\to+\infty.
\label{g.6}
%..............................................................(g.6)
\end{equation}

Now let $U$ be the mapping \eqref{d.7} (for $p=1$)
and $\g=\ext(\eta_1\circ U)$
(remind that $\eta_1$ has been introduced  
in \S\ref{subsection.c.3}, and 
$\ext$ --- in \S\ref{subsection.c.4}).
Note that for any $W\in Y_1$,
$$
\det W=\exp\left(i\int_0^{2\pi} f(e^{i\t})d\t\right),\quad
f\in\pi_1^{-1}(\eta_1(W)).
$$
Thus, it is clear that with the choice \eqref{g.6} of the branch,
one has 
$$
\arg\det M(z;H,H_0)=\int_0^{2\pi}\wt\g(e^{i\t};t)d\t,\quad
z=\l+i(1-t)t^{-1},
$$
where $\wt\g$ is the lift of $\g$
with the initial  condition $\wt\g(0)=0$.
This proves the first of the equalities \eqref{g.3}.
The second one follows from \thm{e.1} after the change of 
variables $t=\cot(\t/2)$.
\qed

%***********************************************************
\section{The invariance principle for $\mu$}
\label{section.f}
%************************************************************

%======================================================
\subsection{Statement of results}
\label{subsection.f.1}
%======================================================

Let $H_0$ and $H$ be self-adjoint operators 
in a Hilbert space $\H$.
Fix $\l\in\R$. 
In this section we prove the invariance principle \eqref{a.13}
for the function $\mu$. 
We find it more natural to prove it in the following form:
\begin{equation}
\mu(\t;f_1(\l),f_1(H),f_1(H_0))=
\mu(\t;f_2(\l),f_2(H),f_2(H_0)),
\quad \t\in(0,2\pi).
\label{f.0a}
%.........................................................(f.0a)
\end{equation}
The functions $f_1$, $f_2$ in \eqref{f.0a} are supposed to 
satisfy \cnd{a.1} (with the same $\Omega\supset\s(H_0)\cup\s(H)$
for $f_1$ and $f_2$ and with $\l$ from \eqref{f.0a}).
\begin{theorem}
\label{t.f.2}
Let $\Omega\subset\R$ be a Borel set, 
$\s(H_0)\cup\s(H)\subset\Omega$, and let the functions 
$f_1$, $f_2$ satisfy \cnd{a.1} with $\l\in\Omega$.
Let the two pairs of operators 
$f_j(H_0),f_j(H)$, $j=1,2$, 
satisfy \cnd{d.1}(i) (for $p=\infty$).
Then:

{\rm (i)} \cnd{d.3} (for $p=\infty$) holds true 
for the pair $f_1(H_0),f_1(H)$ 
at the point $f_1(\l)$ if and only if 
it holds true for the pair $f_2(H_0),f_2(H)$ 
at the point $f_2(\l)$.

{\rm (ii)} If for $j=1,2$ 
\cnd{d.3} (for $p=\infty$) holds true 
for the pair $f_j(H_0),f_j(H)$ 
at the point $f_j(\l)$, then 
\begin{equation}
\xinflim_{y\to0+}\eta_\infty
\bigl(M(f_1(\l)+iy;f_1(H),f_1(H_0))\bigr)
=
\xinflim_{y\to0+}\eta_\infty
\bigl(M(f_2(\l)+iy;f_2(H),f_2(H_0))\bigr).
\label{f.1}
%.........................................................(f.1)
\end{equation}
\end{theorem}

Suppose that under the hypothesis of \thm{f.2}, 
the two pairs of operators $f_j(H_0),f_j(H)$, $j=1,2$, 
satisfy the \cnd{d.1}(ii) (for $p=\infty$).
Then $\mu(\cdot;f_j(\l),f_j(H),f_j(H_0))$ is well defined for 
$j=1,2$. The relation \eqref{f.1} leads to the invariance principle 
\eqref{f.0a} modulo $\Z$. 
In order to obtain the invariance principle in the full scale,
we have to replace \cnd{d.1} by a pair of slightly more restrictive 
conditions.

For $z\in\C$, $z\notin\R_-:=\{z\mid \Im z=0,\Re z<0\}$,
let us fix the branch of $\arg z$, say, by 
\begin{equation}
\arg z\in(-\pi,\pi),\quad
z\in\C\setminus\R_-.
\label{f.2}
%.........................................................(f.2)
\end{equation}
\begin{condition}
\label{condition.f.3}
For a pair of self-adjoint operators $H_0$, $H$, one has:

{\rm (i)}
for any $z\in\C_+$, 
\begin{equation}
\arg(H-zI)-\arg(H_0-zI)\in\SS_\infty(\H);
\label{f.3}
%.........................................................(f.3)
\end{equation}

{\rm (ii)}
for any $\l\in\R$,  
\begin{equation}
\lim_{y\to+\infty}
\norm{\arg(H-(\l+iy)I)-\arg(H_0-(\l+iy)I)}=0.
\label{f.4}
%.........................................................(f.4)
\end{equation}
\end{condition}
\begin{proposition}
\label{p.f.4a}
If for the pair  $H_0$, $H$
\cnd{f.3}(i) holds, then \cnd{d.1}(i) holds.
If \cnd{f.3}(ii) holds, then \cnd{d.1}(ii) holds.
\end{proposition}
\begin{theorem}
\label{t.f.5}
Let $\Omega\subset\R$ be a Borel set, 
$\s(H_0)\cup\s(H)\subset\Omega$, and let the functions 
$f_1$, $f_2$ satisfy \cnd{a.1} with $\l\in\Omega$.
Let, for $j=1,2$, the pair of operators $f_j(H_0),f_j(H)$
satisfy \cnd{f.3} and 
\cnd{d.3} (for $p=\infty$) at the point $f_j(\l)$.
Then the invariance principle \eqref{f.0a} holds.
\end{theorem}
Let us give a sufficient condition for \cnd{f.3}.
\begin{theorem}
\label{t.f.5a}
Let the operators $H_0$, $G$, $J$ be as described in 
\S\ref{subsection.b.2}; assume \eqref{b.2} and let
$H=H(H_0,G,J)$. 
Then \cnd{f.3} holds for the pair $H_0$, $H$.
\end{theorem}

%====================================================
\subsection{Corollaries}
\label{subsection.f.1a}
%====================================================

Theorems \ref{t.e.1}, \ref{t.g.1} and \ref{t.f.5} imply the
following statement, which is the central result of this paper.

\begin{theorem}
\label{t.f.5b}
Let the operators $H_0$, $G$, $J$ be as described in 
\S\ref{subsection.b.2}; assume \eqref{b.2} and let
$H=H(H_0,G,J)$.
Suppose that for an open interval $\d\subset\R$ 
the inclusion \eqref{d.9} holds. 
Further, let $\Omega\subset\R$ be a Borel set, 
$\s(H_0)\cup\s(H)\subset\Omega$, and let a function 
$f$ satisfy \cnd{a.1} for all $\l\in\d$. 
Suppose that 
$$
f(H)-f(H_0)\in\SS_1(\H).
$$
Then for a.e. $\l\in\d$, the representation \eqref{a.10}
holds true.
\end{theorem}
{\em Proof}
First note that the limit $T(\l+i0)$ exists in $\SS_\infty(\K)$ 
by \prp{d.7} and the pair 
$\Xi(\J+A(\l+i0)+tB(\l+i0)),\Xi(\J)$ is Fredholm by 
\thm{e.1}.
Further, by \thm{f.5a}, both the pair $H_0$, $H$, and the pair 
$f(H_0)$, $f(H)$ satisfy \cnd{f.3}.
By \prp{d.7}, the pair $H_0$, $H$ satisfies \cnd{d.3} 
(for $p=\infty$) for a.e. $\l\in\d$ and the pair 
$f(H_0)$, $f(H)$ satisfies \cnd{d.3} 
(for $p=\infty$) for a.e. $\l\in\R$.
Thus, we can apply \thm{f.5}, which yields
$$
\mu(\t;f(\l),f(H),f(H_0))=
\mu(\t;\l,H,H_0),\quad \text{ a.e. }\l\in\d.
$$
By \thm{e.1}, one has
$$
\mu(\theta;\l,H,H_0)=
\iindex\bigl(\Xi(\J), \Xi(\J+A(\l+i0)+\cot(\t/2)B(\l+i0))\bigr).
$$
Applying \thm{g.1} to the pair $f(H_0)$, $f(H)$, we get 
$$
\lim_{y\to0+}\arg\D_{f(H)/f(H_0)}(\l'+iy)=
-\frac12\int_0^{2\pi}\mu(\t;\l',f(H),f(H_0))d\t,
\text{ a.e. }\l'\in\R.
$$
Combining the last three equalities and 
the Krein's formula \eqref{a.3}
and making the change 
of variables $t=\cot(\t/2)$ in the resulting integral,
we get \eqref{a.10}.
\qed

As in \S\ref{subsection.e.1},
for the perturbations of a definite sign 
the representation \eqref{a.10} takes the form
\begin{align}
\xi(f(\l);f(H),f(H_0))&=
\frac{1}{\pi}\int_{-\infty}^\infty\frac{dt}{1+t^2}
\rank E_{A(\l+i0)+tB(\l+i0)}((-\infty,-1)),\quad &J=I,
\label{f.5b}
%........................................................(f.5b)
\\
\xi(f(\l);f(H),f(H_0))&=
-\frac{1}{\pi}\int_{-\infty}^\infty\frac{dt}{1+t^2}
\rank E_{A(\l+i0)+tB(\l+i0)}([1,\infty)),
\quad &J=-I.
\label{f.5c}
%........................................................(f.5c)
\end{align}

The representations \eqref{f.5b}, \eqref{f.5c}
have been originally proven in \cite{Push}
in the following particular case. It was assumed that the operator 
$H_0$ is semibounded from below and $f(\l)=(\l-a)^{-l}$, $l>0$,
$a<\inf(\s(H_0)\cup\s(H))$. Instead of  \eqref{d.9},
it was supposed  that $G(H_0-aI)^{-m}\in\SS_2$ for some $m>0$.
The proof was heavily based upon the particular form of 
the function $f$ and used the results of \cite{Koplienko}.

Note that the SSF is non-negative in \eqref{f.5b} and non-positive
in \eqref{f.5c}. This fact itself is already non-trivial.
In the case $f(\l)=\l$, it has been proven by M.~G.~Krein
in the original paper \cite{Krein}, but very few generalisations
for $f(\l)\not=\l$ have been known so far
(see \cite[\S8.10]{Yafaev} for the discussion).

%========================================================
\subsection{Auxiliary statements}
\label{subsection.f.2}
%========================================================
\begin{lemma}
\label{l.f.6}
Let $M_j=M_j^*\in{\mathcal B}(\H)$, $j=0,1$.
Then, for any $t\in\R$, 

{\rm (i)} one has
$$
\norm{e^{itM}-e^{itM_0}}\leq\abs{t}\norm{M-M_0};
$$

{\rm (ii)} if $M-M_0\in\SS_\infty$, then 
$e^{itM}-e^{itM_0}\in\SS_\infty$.
\end{lemma}
{\it Proof }
Immediately follows from the representation
$$
e^{itM}-e^{itM_0}=
ie^{itM}\int_0^t e^{-isM}(M-M_0)e^{isM_0}ds.\quad
\qed
$$

Recall that we have fixed 
the branch of the argument by \eqref{f.2}.
\begin{lemma}
\label{l.f.7}
Let the functions $f_1$, $f_2$ satisfy \cnd{a.1}
at a point $\l\in\Omega$.
Then, for any self-adjoint operator $H$
such that $\s(H)\subset\Omega$, one has
\begin{equation}
\lim_{y\to0+}\norm{
\arg(f_2(H)-f_2(\l)I-iyf'_2(\l)I)-
\arg(f_1(H)-f_1(\l)I-iyf'_1(\l)I)}=0.
\label{f.6}
%.........................................................(f.6)
\end{equation}
\end{lemma}
{\it Proof }
1.
First let us denote $g_j(x)=(f_j(x)-f_j(\l))/f'_j(\l)$,
$j=1,2$ and without loss of generality assume that $\l=0$.
Clearly, we get $g_j(0)=0$, $g'_j(0)=1$, $j=1,2$, 
and we have to prove that 
$$
\lim_{y\to0+}\norm{\arg(g_2(H)-iyI)-\arg(g_1(H)-iyI)}=0,
$$
which reduces to
\begin{equation}
\lim_{y\to0+}\sup_{x\in\R}
\abs{\arg(g_2(x)-iy)-\arg(g_1(x)-iy)}=0.
\label{f.7}
%.........................................................(f.7)
\end{equation}
It is sufficient to prove the following two relations:
\begin{gather}
\lim_{y\to0+}\sup_{\abs{x}>\d}
\abs{\arg(g_2(x)-iy)-\arg(g_1(x)-iy)}=0
\quad\text{ for any }\d>0,
\label{f.9}
%.........................................................(f.9)
\\
\lim_{x\to0}\sup_{y>0}
\abs{\arg(g_2(x)-iy)-\arg(g_1(x)-iy)}=0.
\label{f.8}
%.........................................................(f.8)
\end{gather}

2.
Let us prove \eqref{f.9}. Clearly, by
\cnd{a.1}(ii), one has 
$\sup_{\abs{x}>\d}(1/\abs{g_j(x)})<\infty$, $j=1,2$.
Thus, as $y\to0+$,
$$
\arg(g_2(x)-iy)-\arg(g_1(x)-iy)=\arg\br{1-(iy/g_2(x))}-
\arg\br{1-(iy/g_1(x))}=O(y)
$$
uniformly in $\abs{x}>\d$.

3.
Let us prove \eqref{f.8}. By \cnd{a.1}(i), one has for $x\to0$:
\begin{align*}
\arg(g_2(x)-iy)-\arg(g_1(x)-iy)&=
\arg(x+o(x)-iy)-\arg(x+o(x)-iy)\\
&=\arg\br{1-i(y/x)+o(1)}-\arg\br{1-i(y/x)+o(1)}=o(1)
\end{align*}
uniformly in $y>0$.
\qed

%===============================================================
\subsection{Proof of \prp{f.4a} and
Theorems \protect\ref{t.f.2}, \ref{t.f.5}}
\label{subsection.f.3}
%===============================================================
{\bf 1.} {\it Proof of \prp{f.4a} }
First note that 
$$
M(z;H,H_0)-I=\exp(-2i\arg(H-zI))
\br{\exp(2i\arg(H_0-zI))-\exp(2i\arg(H-zI))}.
$$
Thus, by \lma{f.6}(ii), 
\eqref{f.3} implies  \eqref{d.3} (with $p=\infty$).
The inclusion \eqref{d.3} is equivalent to \eqref{d.8aa}.
Similarly, by \lma{f.6}(i), \eqref{f.4} implies 
\eqref{d.4a} (with $p=\infty$), and \eqref{d.4a}
is equivalent to \eqref{d.3a}.
\qed

{\bf 2.}
{\it Proof of \thm{f.2} }
As in the proof of \lma{f.7}, we can reduce the problem 
to the case $\l=0$, $f_j(0)=0$, $f'_j(0)=1$, $j=1,2$.
Further, for $x\in\R$ and $y>0$ denote 
$$
A(x;y)=\frac{(f_2(x)+iy)(f_1(x)-iy)}{(f_2(x)-iy)(f_1(x)+iy)}=
\exp\br{2i\arg(f_1(x)-iy)-2i\arg(f_2(x)-iy)}.
$$
One has 
$$
M(iy;f_2(H),f_2(H_0))=
A(H;y)M(iy;f_1(H),f_1(H_0))
(A(H_0;y))^*.
$$
By \lma{f.7} and \lma{f.6}(i),
$$
\lim_{y\to0+}\norm{A(H;y)-I}=
\lim_{y\to0+}\norm{A(H_0;y)-I}=0.
$$
Therefore,
$$
\lim_{y\to0+}\norm{M(iy;f_2(H),f_2(H_0))-
M(iy;f_1(H),f_1(H_0))}=0.
$$
By \lma{c.5}, this proves the theorem.
\qed

{\bf 3.}
{\it Proof of \thm{f.5} }
1. For $j=1,2$, let 
$U_j$ be the mapping \eqref{d.7} (for $p=\infty$),
corresponding to the pair of operators $f_j(H_0),f_j(H)$ and 
the spectral parameter $f_j(\l)$.
Let $\g_j=\ext(\eta_\infty\circ U_j)$
(recall that $\eta_\infty$ has been introduced  
in \S\ref{subsection.c.3}, and 
$\ext$ --- in \S\ref{subsection.c.4}).
Clearly, $\g_1(0)=\g_2(0)$.
By \thm{f.2}, ${\g_1}(1)={\g_2}(1)$.
Below we explicitely construct a homotopy between $\g_1$
and $\g_2$.
By \prp{c.5}, the existence of a homotopy between $\g_1$ and $\g_2$ 
implies that
$$
\spflow(z;U_1)=\spflow(z;U_2),\quad z\in\T\setminus\{1\},
$$
and \eqref{f.0a} follows.

2. 
As in the proof of \lma{f.7}, we reduce the problem 
to the case when $\l=0$, $f_j(0)=0$, $f'_j(0)=1$, $j=1,2$.
Further, for $x\in\R$ and $t\in(0,1)$ denote 
$$
h_j(x;t):=\arg(f_j(x)-i(1-t)t^{-1}),\quad j=1,2.
$$
For $x\in\R$, $s\in[0,1]$ and $t\in(0,1)$ denote
\begin{align*}
A(x;t,s)&:=\exp\bigl(2is(h_1(x;t)-h_2(x;t))\bigr),\\
M(t,s)&:=A(H;t,s)M(i(1-t)t^{-1};f_1(H),f_1(H_0))
(A(H_0;t,s))^*.
\end{align*}
It is straightforward to see that 
\begin{equation}
\begin{split}
M(t,0)=M(i(1-t)t^{-1};f_1(H),f_1(H_0)),\\
M(t,1)=M(i(1-t)t^{-1};f_2(H),f_2(H_0)).
\end{split}
\label{f.12}
%.........................................................(f.12)
\end{equation}

3. Let us check that $M(t,s)-I\in \SS_\infty(\H)$ for all
$(t,s)\in(0,1)\times[0,1]$.
By  \cnd{f.3}(i), one has 
$h_j(H;t)-h_j(H_0;t)\in\SS_\infty(\H)$
for all $t\in(0,1)$ and $j=1,2$.
By \lma{f.6}(ii), this implies that 
$$
\exp(2ish_j(H;t))-\exp(2ish_j(H_0;t))\in\SS_\infty(\H),
\quad(t,s)\in(0,1)\times[0,1],\quad j=1,2
$$
and therefore
$$
A(H;t,s)-A(H_0;t,s)\in\SS_\infty(\H),
\quad(t,s)\in(0,1)\times[0,1].
$$
From here it is easy to infer that $M(t,s)-I\in \SS_\infty(\H)$.

4.
Define the mapping $\G:[0,1]\times[0,1]\to X_\infty$ by 
\begin{align*}
\G(t,s)&=\eta_\infty(M(t,s)),\quad t\not=0,1;\\
\G(0,s)&=0;\\
\G(1,s)&={\g_1}(1)(=\g_2(1)).
\end{align*}
Let us prove that $\G$ is a homotopy 
between ${\g_1}$ and ${\g_2}$.
By \eqref{f.12}, 
$\G(t,0)={\g_1}(t)$ 
and
$\G(t,1)={\g_2}(t)$
for all $t\in[0,1]$. 
It remains to check that the mapping $\G$ is continuous.

5. 
First let us check that the mapping
$$
(0,1)\times[0,1]\ni(t,s)\mapsto M(t,s)-I\in \SS_\infty(\H)
$$
is continuous. 
By \prp{d.1a}(iii), $M(i(1-t)t^{-1};f_1(H),f_1(H_0))$
depends continuously on $t\in(0,1)$ in the operator norm.
It can also be checked explicitly that the mapping 
$$
(0,1)\times[0,1]\ni(t,s)\mapsto A(\cdot;t,s)\in C(\R)
$$
is continuous and therefore $A(H;t,s)$ and $A(H_0;t,s)$ 
depend continuously on $(t,s)$ 
in the operator norm.

6.
Let us check the continuity of $\G$ at $t=0$.
Let us prove that 
$$
\lim_{t\to0+}\sup_{s\in[0,1]}\norm{M(t,s)-I}=0.
$$
\cnd{f.3}(ii) implies that 
$$
\lim_{y\to+\infty}\norm{M(iy;f_1(H),f_1(H_0))-I}=0,
$$ 
and therefore it suffices to prove that 
$$
\lim_{t\to0+}\sup_{s\in[0,1]}
\norm{A(H;t,s)-A(H_0;t,s)}=0,
\quad j=1,2.
$$
By \lma{f.6}(i),
the last relation follows again from \cnd{f.3}(ii).

7.
Let us check the continuity of $\G$ at $t=1$.
Let us prove that 
\begin{equation}
\lim_{t\to1-}\sup_{s\in[0,1]}\norm{M(t,s)-M(t,0)}=0.
\label{f.13}
%.........................................................(f.13)
\end{equation}
It follows from \lma{f.7} and \lma{f.6}(i) that 
$$
\lim_{t\to1-}\sup_{s\in[0,1]}\norm{A(H_0;t,s)-I}=
\lim_{t\to1-}\sup_{s\in[0,1]}\norm{A(H;t,s)-I}=0.
$$
This implies \eqref{f.13}.
By \lma{c.5}, it follows that $\G$ is continuous 
at $t=1$.
\qed

%===============================================================
\subsection{Proof of \thm{f.5a}}
\label{subsection.f.4}
%===============================================================

\begin{lemma}
\label{l.h.1}
Let $H_0$ be a self-adjoint operator in $\H$ and 
$K$ be a compact operator. Then, for any $r>0$ and 
$\psi\in\H$ one has
\begin{gather}
\int_r^\infty
\norm{K(\abs{H_0}+I)^{1/2}(H_0-itI)^{-1}\psi}^2dt\leq
C_{\ref{h.1}}(r;K)\norm{\psi}^2,
\label{h.1}
%.........................................................(h.1)
\\
\intertext{where}
\lim_{r\to\infty}C_{\ref{h.1}}(r;K)=0.
\label{h.2}
%.........................................................(h.2)
\end{gather}
\end{lemma}
{\it Proof }
1. Below we prove the following two facts:

(i) the relations \eqref{h.1}, \eqref{h.2} hold for any finite 
rank operator $K$;

(ii) for any bounded operator $K$ and any $\psi\in\H$ one has
\begin{equation}
\int_1^\infty\norm{K(\abs{H_0}+I)^{1/2}(H_0-itI)^{-1}\psi}^2dt\leq
C_{\ref{h.3}}\norm{K}^2\norm{\psi}^2,
\label{h.3}
%.........................................................(h.3)
\end{equation}
where $C_{\ref{h.3}}$ is a universal constant.

Approximating a compact operator $K$ by finite rank operators,
it is easy to obtain the assertion of the lemma from (i), (ii).

2. Let us prove (i). Clearly, it is sufficient to consider 
a rank one operator $K=(\cdot,\f)\chi$, 
$\norm{\f}=\norm{\chi}=1$.
Let $d\mu_\f(\l):=d(E_{H_0}((-\infty,\l))\f,\f)$ be 
the spectral measure of $H_0$, associated with the vector $\f$.
One has:
\begin{align*}
\int_r^\infty\norm{K(\abs{H_0}+I)^{1/2}(H_0-itI)^{-1}\psi}^2dt&=
\int_r^\infty\abs{((\abs{H_0}+I)^{1/2}(H_0-itI)^{-1}\psi,\f)}^2dt\\
&\leq\norm{\psi}^2
\int_r^\infty\norm{(\abs{H_0}+I)^{1/2}(H_0+itI)^{-1}\f}^2dt\\
&=\norm{\psi}^2\int_r^\infty dt\int_\R
\frac{\abs{\l}+1}{\l^2+t^2}d\mu_\f(\l)\\
&=\norm{\psi}^2\int_\R F(\l,r)d\mu_\f(\l),
\end{align*}
where
$$
F(\l,r)=(\abs{\l}+1)\int_r^\infty\frac{dt}{\l^2+t^2}=
\frac{\abs{\l}+1}{\abs{\l}}\tan^{-1}(\abs{\l}/r).
$$
Clearly, 
\begin{equation}
C_{\ref{h.4}}:=\sup_{r>1}\sup_{\l\in\R}F(\l,r)<\infty,
\label{h.4}
%.........................................................(h.4)
\end{equation}
and $\lim_{r\to\infty}F(\l,r)=0$ for any $\l\in\R$.
Therefore,
$$
\lim_{r\to\infty}\int_\R F(\l,r)d\mu_\f(\l)=0
$$
and we arrive at \eqref{h.1}, \eqref{h.2} with 
$C_{\ref{h.1}}=\int_\R F(\l,r)d\mu_\f(\l)$.

3. Let us prove (ii). As above, one has:
\begin{align*}
\int_1^\infty\norm{K(\abs{H_0}+I)^{1/2}(H_0-itI)^{-1}\psi}^2dt
&\leq\norm{K}^2
\int_1^\infty\norm{(\abs{H_0}+I)^{1/2}(H_0-itI)^{-1}\psi}^2dt\\
&=\norm{K}^2\int_\R F(\l,1)d\mu_\psi(\l)
\leq C_{\ref{h.4}}\norm{K}^2\norm{\psi}^2,
\end{align*}
and we get \eqref{h.3} with $C_{\ref{h.3}}=C_{\ref{h.4}}$.
\qed

{\it Proof of \thm{f.5a} }
1. First of all, note that the conditions \eqref{b.2}
are invariant under the linear transformations $H_0\mapsto aH_0+bI$,
$a,b\in\R$. 
Thus, it is sufficient to prove \eqref{f.3} with $z=i$ and 
\eqref{f.4} --- with $\l=0$.

Next, we will use the integral representation 
$$
\arg(x-iy)=-(\pi/2)+\Im\int_y^\infty\frac{x}{x-it}\frac{dt}{t},
\quad x\in\R,\quad y>0.
$$
In view of this representation, it is sufficient to prove that 
under the assumptions \eqref{b.2}, one has
\begin{gather}
\int_1^R[(H-itI)^{-1}-(H_0-itI)^{-1}]dt\in\SS_\infty(\H)
\quad\text{ for any }R>0,
\label{h.5}
%.........................................................(h.5)
\\
\lim_{r\to\infty}\sup_{R\geq r}
\norm{\int_r^R[(H-itI)^{-1}-(H_0-itI)^{-1}]dt}=0.
\label{h.6}
%.........................................................(h.6)
\end{gather}
By \eqref{b.8}, the inclusion \eqref{d.8aa} (with $p=\infty$)
holds for all $z\in\r(H_0)\cap\r(H)$.
From here we get \eqref{h.5}. 
Thus, it remains to prove \eqref{h.6}.

2.
Let us prove \eqref{h.6}. First, for brevity we denote 
$K:=G(\abs{H_0}+I)^{-1/2}$. Using \eqref{b.8} and 
\lma{h.1}, we obtain the following estimate for any
$\psi,\f\in\H$:
\begin{align*}
\Bigl\lvert\int_r^R&[((H-itI)^{-1}\f,\psi)-((H_0-itI)^{-1}\f,\psi)]dt
\Bigr\rvert\\
&\leq\int_r^R\norm{(\J+T(it))^{-1}}
\norm{G(H_0-itI)^{-1}\f}
\norm{G(H_0-itI)^{-1}\psi}dt\\
&\leq\sup_{t\geq r}\norm{(\J+T(it))^{-1}}
\left(\int_r^\infty
\norm{K(\abs{H_0}+I)^{1/2}(H_0-itI)^{-1}\f}^2dt
\right)^{1/2}\\
&\quad\times
\left(\int_r^\infty
\norm{K(\abs{H_0}+I)^{1/2}(H_0-itI)^{-1}\psi}^2dt
\right)^{1/2}\\
&\leq\sup_{t\geq r}\norm{(\J+T(it))^{-1}}\norm{\f}\norm{\psi}
C_{\ref{h.1}}(r,K),
\end{align*}
which, by \eqref{h.2} and \eqref{d.13b} (with $p=\infty$),
proves \eqref{h.6}.
\qed

%****************************************************************
\section{Appendix: additional properties of the function $\mu$}
\label{section.i}
%****************************************************************

Here we prove formula \eqref{a.11} and explain the relation 
of the function $\mu(\cdot;\l,H,H_0)$ to the eigenvalue counting functions 
of the operators $H_0$ $H$. These results have not been used above 
and are given only in order to clarify the  links of the function $\mu$
to the standard objects of the spectral theory of perturbations.

%===========================================================
\subsection{The function $\mu$ and the spectrum of 
the scattering matrix}
\label{subsection.i.1}
%===========================================================

Let the operators $H_0$, $G$, $J$ be as described in 
\S\ref{subsection.b.2}; assume \eqref{b.2}
and  let $H=H(H_0,G,J)$.
Fix an interval $\D$ in the absolutely continuous spectrum of
$H_0$.
Below we give a criterion for existence of the scattering matrix 
$\calS(\l;H,H_0)$ for a.e. $\l\in\D$,
which can be found, e.g., in 
{\cite[\S5.8]{Yafaev}}.
For technical reasons, we suppose that 
$\Ker G=\{0\}$; this will simplify the statement below.
\begin{proposition}
\label{p.i.4}
Suppose that for  a.e. $\l\in\D$, the limit 
$$
\blim_{y\to0+}T(\l+iy)
$$
exists and $0\in\r(\J+T(\l+i0))$. 
Then the local wave operators $W_\pm(H,H_0;\D)$
exist and are complete. 
For a.e. $\l\in\D$, the scattering matrix 
$\calS(\l;H,H_0)$ is given by
\begin{equation}
\calS(\l;H,H_0)=I-2\pi iZ(\l;G)(\J+T(\l+i0))^{-1}Z^*(\l;G),
\label{i.1}
%..............................................................(i.1)
\end{equation}
where the operator $Z(\l;G)$ satisfies the relation
$$
\pi Z^*(\l;G)Z(\l;G)=B(\l+i0).
$$
\end{proposition}
In this situation, clearly, $\calS(\l;H,H_0)-I\in \SS_\infty$.
Note that under the hypothesis of \prp{i.4}, the 
Assumptions \ref{condition.d.1} and \ref{condition.d.3}
hold for $p=\infty$ and a.e. $\l\in\D$.

Further, by \eqref{d.18a} and \thm{d.10}, one has
$$
\eta_\infty(\calS(\l;H,H_0))=
\eta_\infty(S(\l+i0;H_0,G,J))=
\lim_{y\to0+}\eta_\infty(M(\l+iy;H,H_0)).
$$
Thus, we see that under the hypothesis of \prp{i.4},
for a.e. $\l\in\D$ the relation \eqref{a.11} holds true.

%===========================================================
\subsection{The function $\mu$ on the discrete spectrum}
\label{subsection.i.2}
%===========================================================

{\bf 1.} Let $H_0$, $H$ be self-adjoint operators in $\H$, 
satisfying \cnd{d.1} (with $p=\infty$). 
If $\l\in\R\setminus(\s(H_0)\cup\s(H))$,
then, obviously, \cnd{d.3} is fulfilled and 
$M(\l;H,H_0)=I$.
Therefore, $\mu(\t;\l,H,H_0)$ equals to an integer constant.
Below we discuss the relation of this constant 
to the eigenvalue counting functions of $H_0$ and $H$. 
First we need notation, similar to \eqref{c.4a}, but for 
self-adjoint operators. 
For $\l_1,\l_2\in\R$ and $H=H^*$ we put
\begin{equation*}
N(\l_1,\l_2;H)=\left\{
\begin{array}{lc}
\rank E_H([\l_1,\l_2)),&\l_2>\l_1,\\
0,&\l_2=\l_1,\\ 
-\rank E_H
([\l_2,\l_1)),&\l_1>\l_2.
\end{array}
\right.
\end{equation*}
Remind that \cnd{d.1} implies that $\s_{ess}(H)=\s_{ess}(H_0)$.

\begin{theorem}
\label{t.i.1}
Let $[\l_1,\l_2]\cap\s_{ess}(H_0)=\emptyset$ and 
$\{\l_1,\l_2\}\subset\r(H)\cap\r(H_0)$.
Then, for all $\t\in(0,2\pi)$,
\begin{equation}
\mu(\t;\l_2,H,H_0)-\mu(\t;\l_1,H,H_0)=
N(\l_1,\l_2;H)-N(\l_1,\l_2;H_0).
\label{i.1aa}
%............................................................(i.1aa)
\end{equation}
\end{theorem}

{\bf 2.} Let 
$$
H: [0,1]\ni\a\mapsto H(\a)
$$ 
be a family of self-adjoint operators in $\H$, which satisfies 
the following assumptions:
\begin{gather}
(H(\a)-zI)^{-1}-(H(0)-zI)^{-1}\in\SS_\infty(\H),
\quad \forall z\in\C_+,\quad\a\in[0,1],
\label{i.1a}
%............................................................(i.1a)
\\
\text{the mapping }
[0,1]\ni\a\mapsto (H(\a)-zI)^{-1}\in\calB(\H)
\text{ is continuous for all }z\in\C_+,
\label{i.1b}
%............................................................(i.1b)
\\
\lim_{y\to+\infty}\sup_{\a\in[0,1]}
y\norm{(H(\a)-(\l+iy)I)^{-1}-(H(0)-(\l+iy)I)^{-1}}=0,
\quad\forall\l\in\R.
\label{i.1c}
%............................................................(i.1c)
\end{gather}
By \eqref{i.1a}, the essential spectra of all the operators 
$H(\a)$ coincide.
Suppose that $\D\subset\R\setminus\s_{ess}(H(\a))$.
Below we explain that for $\l\in\D$ the function 
$\mu(\t;\l,H(1),H(0))$ can be considered as the spectral 
flow of the family $H$ through the point $\l$.

In order to define the spectral flow of the family $H$,
let us repeat (without proofs) the basic steps of the construction of 
\S\ref{section.c}.
First let us fix a function space $\wt X$ where the function 
$\spflow(\l;H)$, $\l\in\D$, will belong to. 
Let $\wt X$ be the set of left continuous bounded non-decreasing 
functions $f:\D\to\Z$. 
There is a lot of freedom in choosing the topology in $\wt X$;
let us consider $\wt X$ with the topology, say, induced by 
the embedding $\wt X\subset L_1(\D)$ 
(we could instead take $L_p(\D)$ with any $p<\infty$).
Consider the equivalence relation
$$ 
f\sim g\Longleftrightarrow \exists n\in\Z: 
\forall x\in\D, 
\quad f(x)=g(x)+n. 
$$ 
Let $X$ be the quotient space $\wt X/\!\!\sim$, and let
$\pi:\wt X\to X$ be the corresponding projection.
In the natural way one defines a topology in $X$ and checks that 
$\pi:\wt X\to X$ is a covering.

Further, note that for every $\a\in[0,1]$ and 
$\l_0\in\D\cap\r(H(\a))$, the function 
$N(\l_0,\cdot;H(\a))$ belongs to $\wt X$.
Define the mapping $\g:[0,1]\to X$ by 
$$
\g(\a)=\pi\bigl(N(\l_0,\cdot;H(\a))\bigr),
\quad\l_0\in\D\cap\r(H(\a)).
$$
This definition does not depend on the choice of $\l_0$.
Since all the eigenvalues of $H(\a)$ depend continuously on $\a$, 
it follows that $\g$ is continuous.
Let $\wt \g$ be a lift of $\g$ to $\wt X$.
Then we put
\begin{equation}
\spflow(\l;H):=\wt\g(\l;1)-\wt\g(\l;0),\quad \l\in\D.
\label{i.6}
%............................................................(i.6)
\end{equation}
As in \S\ref{subsection.c.4}(3), it is easy to see that 
\begin{equation}
\begin{split} 
\spflow(\l;H)&=\langle
\text{the number of eigenvalues of $H(\a)$ that
cross $\l$ leftwards}\rangle
\\
&\quad-\langle\text{the number of eigenvalues of
$H(\a)$ that cross $\l$ rightwards}\rangle 
\end{split}
\label{i.7}
%...........................................................(i.7)
\end{equation}
as $\a$ grows from $0$ to $1$,
whenever the r.h.s. is well defined.

It follows from \thm{i.1} that $\spflow(\l;H)$ and
$\mu(\t;\l,H(1),H(0))$ differ by a function (of $\l$), 
which is identically equal to an integer number. 
The following theorem shows that 
this number equals zero.

\begin{theorem}
\label{t.i.3}
The mapping 
\begin{equation}
[0,1]\ni\a\mapsto\mu(\t;\cdot,H(\a),H(0))\in L_1(\D)
\label{i.5}
%............................................................(i.5)
\end{equation}
is continuous.
\end{theorem}
Thus, the mapping \eqref{i.5} is a lift of $\g$ and therefore,
\begin{equation}
\mu(\t;\l,H(1),H(0))=
\spflow(\l;H),
\quad\l\in\D.
\label{i.8}
%...........................................................(i.8)
\end{equation}

As a typical example, consider the family 
$H(\a)=H(H_0,\sqrt{\a}G,J)$, where the operators
$H_0$, $G$, $J$ satisfy \eqref{b.2}. 
It is easy to see that in this case the assumptions 
\eqref{i.1a}--\eqref{i.1c} hold.
Moreover, the eigenvalues of $H(\a)$ in the gaps 
depend analytically on $\a$, and therefore 
the r.h.s. of \eqref{i.7}
is well defined (see, e.g., \cite{Safronov}).

%===========================================================
\subsection{Proofs of Theorems \ref{t.i.1}, \ref{t.i.3}}
\label{subsection.i.3}
%===========================================================

{\it Proof of \thm{i.1} }
1. Let us first prove that if $[\l_1,\l_2]\subset\r(H_0)\cap\r(H)$,
then 
\begin{equation}
\mu(\t;\l_1,H,H_0)=\mu(\t;\l_2,H,H_0).
\label{i.2}
%............................................................(i.2)
\end{equation}
For $j=1,2$, let $\g_j:[0,1]\to X_\infty$ be the mapping 
\begin{align*}
\g_j(0)&=0,\\
\g_j(t)&=\eta_\infty(M(\l_j+i(1-t)t^{-1};H,H_0)),\quad
t\in (0,1].
\end{align*}
We need to check that $\g_1$ and $\g_2$ are homotopic.
Define the mapping $\G:[0,1]\times[\l_1,\l_2]\to X_\infty$ by 
\begin{align*}
\G(0,\l)&=0,\quad \l\in[\l_1,\l_2];\\
\G(t,\l)&=\eta_\infty(M(\l+i(1-t)t^{-1};H,H_0)),\quad
(t,\l)\in(0,1]\times[\l_1,\l_2].
\end{align*}
Similarly to the proof of \thm{f.5}, one easily checks that 
$\G$ is a homotopy between $\g_1$ and $\g_2$.

2. It remains to check that for all $\l\in\R\setminus\s_{ess}(H_0)$,
one has 
\begin{equation}
\mu(\t;\l+0,H,H_0)-\mu(\t;\l-0,H,H_0)=
\rank E_H(\{\l\})-\rank E_{H_0}(\{\l\}).
\label{i.3}
%............................................................(i.3)
\end{equation}
Without the loss of generality assume that $\l=0$.
Choose $\e>0$ small enough so that there is no spectrum 
of $H$ and $H_0$ in $[-\e,0)\cup(0,\e]$.
We are going to prove that 
$$
\mu(\t;\e,H,H_0)-\mu(\t;-\e,H,H_0)=
\rank E_H(\{0\})-\rank E_{H_0}(\{0\}).
$$
In order to do this, consider the path $\b_1:[0,\pi]\to X_\infty$,
$$
\b_1(t)=\eta_\infty(M(\e e^{it};H,H_0)).
$$
Clearly, $\b_1(0)=\b_1(\pi)=0$. 
Further, consider the paths $\g_\pm:[0,1]\to X_\infty$,
\begin{align*}
\g_\pm(0)&=0,\\
\g_\pm(t)&=\eta_\infty(M(\pm\e+i(1-t)t^{-1};H,H_0)), \quad
t\in (0,1].
\end{align*}
It is easy to see that the catenation $\g_+\cdot\b_1$ 
is homotopic to $\g_-$. Therefore, it is sufficient to prove 
that 
\begin{equation}
\wt\b_1(\t;\pi)-\wt\b_1(\t;0)
=\rank E_{H_0}(\{0\})-\rank E_{H}(\{0\}),
\label{i.4}
%............................................................(i.4)
\end{equation}
where $\wt\b_1$ is a lift of $\b_1$.

3. In order to prove \eqref{i.4}, we are going to check that 
$\b_1$ is homotopic to the following path 
$\b_2:[0,\pi]\to X_\infty$:
$$
\b_2(t):=\eta_\infty\bigl(
(E_{H}(\R\setminus\{0\})+e^{-2it}E_{H}(\{0\}))
(E_{H_0}(\R\setminus\{0\})+e^{2it}E_{H_0}(\{0\}))
\bigr).
$$
It is clear that for a lift $\wt\b_2$ of $\b_2$, one has 
$$
\wt\b_2(\t;\pi)-\wt\b_2(\t;\pi)=
\rank E_{H_0}(\{0\})-\rank E_{H}(\{0\}),
$$
which implies \eqref{i.4}.

4. The homotopy $\G:[0,\pi]\times[0,1]\to X_\infty$ between
$\b_1$ and $\b_2$ is given by 
\begin{align*}
\G(t,s)&=\eta_\infty(U(t,s)),\\
U(t,s)&=\left(\frac{H-s\e e^{-it}I}{H-s\e e^{it}I}
E_{H}(\R\setminus\{0\})+e^{-2it}E_{H}(\{0\})\right)\\
&\quad\times
\left(\frac{H_0-s\e e^{it}I}{H_0-s\e e^{-it}I}
E_{H_0}(\R\setminus\{0\})+e^{2it}E_{H_0}(\{0\})\right)
\qed
\end{align*}

{\it Proof of \thm{i.3} }
1.
First let us prove the following statement.
Fix $\l\in\r(H_0)$ and consider 
$\mu(\t;\l,H(\a),H(0))$ as the function of $\a$. 
Let $\d=[\a_1,\a_2]$ be an interval such that 
$\l\in\r(H(\a))$ for all $\a\in\d$.
Then 
$$
\mu(\t;\l,H(\a_1),H(0))=
\mu(\t;\l,H(\a_2),H(0)).
$$

For $j=1,2$ let $U_j$ be the mapping
\eqref{d.7} (with $p=\infty$)
for the pair $H(0)$, $H(\a_j)$, and let 
$\g=\ext(\eta_\infty\circ U)$.
We need to prove that $\g_1$ and $\g_2$ are homotopic.
Using \eqref{i.1a}--\eqref{i.1c}, one easily checks that 
the mapping $\G:[0,1]\times\d\to X_\infty$, given by 
\begin{align*}
\G(0,\a)&=0,\\
\G(t,\a)&=\eta_\infty(M(\l+i(1-t)t^{-1};H(\a),H(0))),
\end{align*}
is a homotopy between $\g_1$ and $\g_2$.

2.
Fix $\a_0\in[0,1]$; let the neighbourhood $\o\subset[0,1]$ 
of $\a_0$ be small enough so that there exists $\l_0\in\D$, 
$\l_0\in\r(H(\a))$ for all $\a\in\o$.
As we have seen above, one has
$$
\mu(\t;\l_0,H(\a),H(0))=
\mu(\t;\l_0,H(\a_0),H(0)), \quad \a\in\omega.
$$
Therefore, by \thm{i.1},
\begin{multline*}
\mu(\t;\l,H(\a),H(0))-\mu(\t;\l,H(\a_0),H(0))=
(\mu(\t;\l,H(\a),H(0))-\mu(\t;\l_0,H(\a),H(0)))\\
-(\mu(\t;\l,H(\a_0),H(0))-\mu(\t;\l_0,H(\a_0),H(0)))
=N(\l_0,\l;H(\a))-N(\l_0,\l;H(\a_0)).
\end{multline*}
Since there are only finitely many eigenvalues of $H(\a)$ in $\D$ 
and they depend continuously on $t$, we conclude that 
$$
\lim_{\a\to\a_0}
\norm{N(\l_0,\cdot;H(\a))-N(\l_0,\cdot;H(\a_0))}_{L_1(\D)}=0.
$$
This implies \eqref{i.5}.
\qed

\section*{Acknowledgement}
The author is grateful to K.~A. Makarov for useful discussions.

\end{document}